\documentclass{article}
\usepackage{amssymb}
\usepackage{amsmath}
\usepackage{amsxtra}
\usepackage{mathtools}
\usepackage[all,arrow]{xy}

\usepackage{geometry}
\usepackage{threeparttable}
\usepackage{diagbox}
\usepackage{paralist,bbding,pifont}

\usepackage[thmmarks]{ntheorem}
{
	\theoremstyle{nonumberplain}
	\theoremheaderfont{\bfseries}
	\theorembodyfont{\normalfont}
	\theoremsymbol{\mbox{$\Box$}}
	\newtheorem{pf}{Proof}
}
{
	\theoremstyle{numberplain}
	\theoremheaderfont{\bfseries}
	\theorembodyfont{\normalfont}
	\newtheorem{rmk}{Remark}
}
{
	\theoremstyle{numberplain}
	\theoremheaderfont{\bfseries}
	\theorembodyfont{\normalfont}
	\newtheorem{eg}{Example}
}
{
	\theoremstyle{numberplain}
	\theoremheaderfont{\bfseries}
	\theorembodyfont{\normalfont}
	\newtheorem*{ackn}{Acknowledgements.}
}
\newtheorem{thm}{Theorem}

\newtheorem{prop}{Proposition}
\newtheorem{lem}{Lemma}
\newtheorem{Q}{Problem}

\usepackage{hyperref}
\hypersetup{hypertex=true,
	colorlinks=true,
	linkcolor=blue,
	anchorcolor=blue,
	citecolor=blue}

\title{Free circle actions on certain simply connected $7-$manifolds}
\author{Fupeng Xu}
\date{}

\begin{document}
	
\maketitle
\begin{abstract}
	In this paper, we determine for which nonnegative integers $k$, $l$ and for which homotopy $7-$sphere $\Sigma$ the manifold $kS^{2}\times S^{5}\#lS^{3}\times S^{4}\#\Sigma$ admits a free smooth circle action. 
\end{abstract}

%\tableofcontents

\section{Introduction}
In this paper, all manifolds are closed, connected and smooth, and group actions on manifolds are smooth if no otherwise is stated. The standard unit $n-$sphere in the Euclidean space $\mathbb{R}^{n+1}$ is denoted by $S^{n}$. A homotopy $n-$sphere is an $n-$manifold that is homotopy equivalent to $S^{n}$. Given a simply connected $n-$manifold $P$, a \emph{cohomology $P$} is a simply connected $n-$manifold whose integral cohomology ring is isomorphic to that of $P$. Throughout this paper we use (co)homology with integral coefficient unless explicitly stated otherwise.

In topology we ask the following basic problem:
\begin{Q}\label{Problem: main}
	Given a manifold $M$, does $M$ admit a free circle action?
\end{Q}
The solution to this problem and the related topic have applications in the areas of topology and geometry. See \cite[Theorems 7.2.4 and 7.2.5]{Geiges_2008} for application to contact geometry and \cite[Chapeter 6]{Tuschmann2015ModuliSO} for application to differential geometry.

When $M$ is a homotopy sphere of small odd dimensions Problem \ref{Problem: main} was solved by Montgomery and Yang (\cite{MontgomeryYang68}, dimension $7$) and Brumfiel (\cite{Brumfiel69}, dimension $9$, $11$ and $13$). 
When $M$ has certain connectivity and its dimension does not exceed $7$ Problem \ref{Problem: main} was solved by Duan and Liang (\cite{Duan2005CircleBO}, simply connected $5-$manifold), Duan (\cite{duan2022circle}, simply connected $6-$manifold) and Jiang (\cite{RegCirAct2cnt7mfd}, $2-$connected $7-$manifold). 
When $M$ is the connected sum of $S^{q_{i}}\times S^{n-q_{i}}$ or $S^{2}\widetilde{\times}S^{n-2}$, $n\geqslant5$, $2\leqslant q_{i}\leqslant n-2$ and its Betti numbers satisfy some extra conditions, Galaz-Garc\'{i}a and Reiser showed that the answer to Problem \ref{Problem: main} is affirmative (\cite[Theorem E]{galazgarcía2023free}).

%In \cite{MontgomeryYang68} Montgomery and Yang determined which homotopy $7-$spheres admits free smooth circle action. In \cite{Duan2005CircleBO} Duan and Liang classified all simply connected $5-$dimensional smooth manifolds admitting a free smooth circle action. Jiang determined homeomorphism (resp. diffeomrphism) type of $2-$connected $7-$dimensional manifolds admitting regular (resp. smooth free) circle action in \cite{RegCirAct2cnt7mfd}. And Duan solved the classification of simply connected $6-$dimensional that admits a free smooth action in \cite{duan2022circle}. As for higher dimensions, the result of \cite{galazgarcía2023free}, Theorem E tells us the existence of free smooth circle action on a certain type of manifolds that can be writted as connected sum of products of two spheres with dimension both greater than $1$, not restricting to one particular dimension, though only standard smooth structures are concerned and no exotic structure (in particular, no connected sum with exotic spheres) occurs.\\

In this paper we answer Problem \ref{Problem: main} for those simply connected $7-$manifolds that can be written as connected sums of $S^{2}\times S^{5}$, $S^{3}\times S^{4}$ and homotopy $7-$spheres. The connected sum of two $7-$manifolds $M_{1}$ and $M_{2}$ is denoted by $M_{1}\#M_{2}$, and the connected sum of $k$ copies of $7-$manifold $M$ is denoted by $kM$. Our main result is as follows.
\begin{thm}\label{main}
	Let $k,\ l\in\mathbb{N}$ and let $\Sigma$ be a homotopy $7-$sphere. The manifold $kS^{2}\times S^{5}\#lS^{3}\times S^{4}\#\Sigma$ admits a free circle action if and only if one of the following holds:
	\begin{compactenum}
		\item $k\geqslant 2$;
		\item $k=1$ and $l$ is even;
		\item $k=1$, $l$ is odd and $\Sigma$ admits a free circle action;
		\item $k=0$, $l$ is even and $\Sigma$ admits a free circle action;
		\item $k=0$, $l$ is odd and $\Sigma=S^{7}$.
	\end{compactenum}
\end{thm}
\begin{rmk}
	\cite[Theorem E]{galazgarcía2023free} shows that $kS^{2}\times S^{5}\#lS^{3}\times S^{4}$ admits a free circle action for any $k,l\in\mathbb{N}$. Theorem \ref{main} extends this result to $kS^{2}\times S^{5}\#lS^{3}\times S^{4}\#\Sigma$ when $\Sigma$ is a homotopy $7-$sphere.
\end{rmk}
\begin{rmk}
	The cases 4 and 5 of Theorem \ref{main} are covered by \cite[Theorem 1.3]{RegCirAct2cnt7mfd}. Hence in this paper we focus on the cases 1, 2 and 3 of Theorem \ref{main}.
\end{rmk}
\begin{rmk}
	%The oriented diffeomorphism type of homotopy $7-$spheres form a group $\Theta_{7}$ and the group multiplication is induced by connected sum. The group $\Theta_{7}$ is finite cyclic and has order $28$ (\cite{KervaireMilnor63}).
	Up to orientation-preseriving diffeomorphism there are exactly $28$ homotopy $7-$spheres (\cite{KervaireMilnor63}). They can be distinguished by the $\mu-$invariant, which is an orientation-preserving diffeomorphism invariant for certain manifolds and takes value in $\mathbb{Q}/\mathbb{Z}$. There is a homotopy $7-$sphere $\Sigma_{1}$ %that represents the generator of $\Theta_{7}$
	such that $\mu\left(\Sigma_{1}\right)=\frac{1}{28}\ \mathrm{mod}\ 1$. For $r\in\mathbb{N}$ we set $\Sigma_{r}=r\Sigma_{1}$, then $\mu\left(\Sigma_{r}\right)=\frac{r}{28}\ \mathrm{mod}\ 1$ (\cite[Section 6]{EKinv}) and $\Sigma_{r}$ gives all $28$ homotopy $7-$spheres when $r=0,\ 1,\ \cdots,\ 27$.
\end{rmk}
\begin{rmk}\label{Rmk: homotopy spheres admitting circle actions}
\cite[Section 3]{MontgomeryYang68} shows that $\Sigma_{r}$ admits a free circle action if and only if $r=0$, $4$, $6$, $8$, $10$, $14$, $18$, $20$ or $24$. Equivalently, a homotopy $7-$sphere admits a free circle action if and only if its $\mu-$invariant is $0$, $\pm\frac{4}{28}$, $\pm\frac{6}{28}$, $\pm\frac{8}{28}$, $\pm\frac{10}{28}$ or $\frac{14}{28}\ \mathrm{mod}\ 1$. 
\end{rmk}

The paper is organized as follows. First in Section \ref{Homeom S2 times S5} we study the case $(k,l)=(1,0)$. Next in Section \ref{Homeom S2 times S5 sum S3 times S4} we study the case $(k,l)=(1,1)$, where the argument is parallel to the argument of Section \ref{Homeom S2 times S5}. Then in Section \ref{More cnt sum} we handle the remaining cases $(k,l)$ and thus complete the proof of Theorem \ref{main}, where we use the suspension operation that is developed in \cite{duan2022circle} and \cite{galazgarcía2023free}. Finally in Section \ref{Miscellaneous} we add some other results.

\section{Circle actions on $S^{2}\times S^{5}\#\Sigma$}\label{Homeom S2 times S5}

In this section we determine for which homotopy $7-$sphere $\Sigma$ the manifold $S^{2}\times S^{5}\#\Sigma$ admits a free circle action. The organization is as follows. First in Section \ref{Subsec: Orbit B 1} we deal with a broader class of manifolds, namely the spin cohomology $S^{2}\times S^{5}$. In particular when $\Sigma$ is a homotopy $7-$sphere the manifold $S^{2}\times S^{5}\#\Sigma$ is a spin cohomology $S^{2}\times S^{5}$. We will study the invariants of the orbit space $N$ of a free circle action on a spin cohomology $S^{2}\times S^{5}$. Next in Section \ref{Subsec: s inv 1} we compute the diffeomorphism invariants of those spin cohomology $S^{2}\times S^{5}$ that admits free circle action. We express the diffeomorphism invariants of such manifold $M$ in terms of the invariants of the orbit space $N$. Finally in Section \ref{Subsec: compare s inv 1} we compare the diffeomorphism invariants of $M$ and $S^{2}\times S^{5}\#\Sigma$. We determine among those spin cohomology $S^{2}\times S^{5}$ which are diffeomorphic to $S^{2}\times S^{5}\#\Sigma$. Then we deduce for which homotopy $7-$sphere $\Sigma$ the manifold $S^{2}\times S^{5}\#\Sigma$ admits a free circle action.

\subsection{The orbit space $N$}\label{Subsec: Orbit B 1}
In this section we study the invariants of the orbit space of a free circle action on $M$, where $M$ is a spin cohomology $S^{2}\times S^{5}$. By \cite[Chapter 21, Theorem 21.10]{JohnLee} the orbit space $N$ is a $6-$manifold and the projection $M\xrightarrow{p}N$ is a circle bundle. 

Before we state our results, we introduce some notations.
\begin{compactenum}
	\item When $P$ is an oriented $6-$manifold, there is a symmetric cubic form on $H^{2}(P)$:
	\begin{align*}
		\mu=\mu_{P}:H^{2}(P)\otimes H^{2}(P)\otimes H^{2}(P)&\to \mathbb{Z},\\
		(x,y,z) & \mapsto \left<x\cup y\cup z,[P]\right>.
	\end{align*}
	\item Given $x\in H^{q}(P)$, its reduction in $H^{q}(P;\mathbb{Z}/2)$ is be denoted by $\overline{x}$.
	\item Assume $P$ is an oriented $n-$manifold such that its cohomology groups are torsion free. For $0<q<n$, let $\left(e_{1},\cdots,e_{r}\right)$ be an ordered basis of $H^{q}(P)\cong\mathbb{Z}^{r}$. Then by Poincar\'e duality and universal coefficient theorem $H^{n-q}(P)\cong\mathbb{Z}^{r}$ admits an ordered basis $\left(e^{1},\cdots,e^{r}\right)$ such that $$\left<e_{i}\cup e^{j},[P]\right>=\begin{dcases*}
		1,&if $i=j$,\\
		0,&otherwise.
	\end{dcases*}$$The ordered basis $\left(e^{1},\cdots,e^{r}\right)$ is called the \emph{dual basis} with respect to $\left(e_{1},\cdots,e_{r}\right)$.
\end{compactenum} 
\begin{lem}\label{Lemma: Base space 1}
	Assume the $6-$manifold $N$ is the orbit of a free circle action on a spin cohomology $S^{2}\times S^{5}$. Then $N$ is simply connected and satisfies the following conditions.
	\begin{compactenum}
		\item Its cohomology groups are torsion-free and the only nonzero Betti numbers are $b_{0}=b_{6}=1,\ b_{2}=b_{4}=2$;
		\item Its second cohomology group has a basis $(e,f)$ such that the matrix $M_{e}=\begin{pmatrix}
			A & B \\ B & C
		\end{pmatrix}$ is invertible, where $A=\mu(e,e,e)$, $B=\mu(e,e,f)$ and $C=\mu(e,f,f)$;
		\item Exactly one of the following is true:
		\begin{compactenum}
			\item $N$ is spin and $p_{1}(N)=(24u+4A)e\spcheck+(24v+4D)f\spcheck$ for some $u,v\in\mathbb{Z}$. Meanwhile $A$ is even and $B$, $C$ are odd.
			\item $N$ is nonspin with $w_{2}(N)=\overline{e}$ and $p_{1}(N)=(48u+A)e\spcheck+(24v+3B+6C+4D)f\spcheck$ for some $u,v\in\mathbb{Z}$. Meanwhile $A$, $C$ are even and $B$ is odd.
		\end{compactenum}Here $D=\mu(f,f,f)$ and $\left(e\spcheck,f\spcheck\right)$ is the dual basis of $H^{4}(N)$ with respect to $(e,f)$.
	\end{compactenum}
\end{lem}
%\begin{rmk}
%	Lemma \ref{Lemma: Base space 1} is originally introduced to characterize the orbit space of a free circle action on $S^{2}\times S^{5}\#\Sigma$. In the proof we only used the fact that $S^{2}\times S^{5}\#\Sigma$ is a spin cohomology $S^{2}\times S^{5}$. Meanwhile the necessary condition we deduced is sufficient to characterize the orbit space of a free circle action on a spin cohomology $S^{2}\times S^{5}$ conversely. Hence we divert to characterize the orbit space of a free circle action on a spin cohomology $S^{2}\times S^{5}$ and state Lemma \ref{Lemma: Base space 1} in this way. Its sufficiency will be justified in Lemma \ref{Lemma: Base space 1, if part}.
%\end{rmk}
\begin{pf}
Let $M$ be a spin cohomology $S^{2}\times S^{5}$ that admits a free circle action and let $N$ be the orbit manifold. We shall justify that $N$ is simply connected and satisfies Conditions 1, 2 and 3 of Lemma \ref{Lemma: Base space 1}.

We begin with the cohomology groups of $N$. The free circle action induces a circle bundle $M\xrightarrow{p}N$. Let $e\in H^{2}(N)$ be its Euler class. The long exact sequence of homotopy groups associated to the circle bundle $M\to N$ shows that $N$ is simply connected and there is a short exact sequence $0\to\mathbb{Z}\to\pi_{2}(N)\to\mathbb{Z}\to0$. Hence $N$ is orientable and $\pi_{2}(N)\cong\mathbb{Z}^{2}$. Then $H_{2}(N)\cong\pi_{2}(N)\cong\mathbb{Z}^{2}$ by Hurewicz theorem, and $H^{2}(N)\cong H^{4}(N)\cong\mathbb{Z}^{2}$ by universal coefficient theorem and Poincar\'e duality. From the Gysin sequence associated to the circle bundle $M\xrightarrow{p}N$ we obtain an exact sequence $H^{1}(N)\xrightarrow{\cup e}H^{3}(N)\xrightarrow{p^{*}}H^{3}(M)$. Hence $H^{3}(N)=0$. In summary, the orbit space $N$ is a simply connected $6-$manifold and its cohomology groups are torsion-free with only nonzero Betti numbers $b_{0}=b_{6}=1$, $b_{2}=b_{4}=2$. This justifies Condition 1 in Lemma \ref{Lemma: Base space 1}.

Next we consider the ring structure of $H^{*}(N)$. From the Gysin sequence associated to the circle bundle $M\xrightarrow{p}N$ we obtain another two exact sequences
\begin{align}
	0\to H^{0}(N)\xrightarrow{\cup e}H^{2}(N)\xrightarrow{p^{*}}H^{2}(M)\to 0,\label{SES: e primitve}\\
	0\to H^{2}(N)\xrightarrow{\cup e}H^{4}(N)\to 0.\label{SES: cup e isom}
\end{align}
From the exact sequence \eqref{SES: e primitve} we see $e\in H^{2}(N)\cong\mathbb{Z}^{2}$ is \emph{primitive}, i.e. the subgroup $\left<e\right>$ of $H^{2}(N)$ generated by $e$ is isomorphic to $\mathbb{Z}$ and $H^{2}(N)/\left<e\right>\cong\mathbb{Z}$. Hence we can extend $e$ to an ordered basis $(e,f)$ of $H$. Next we fix a fundamental class $[N]\in H_{6}(N)\cong\mathbb{Z}$ and let $\left(e\spcheck,f\spcheck\right)$ be the dual basis of $H^{4}(N)$ with respect to $(e,f)$. The cubic form $\mu=\mu_{N}$ is determined by the integers $A=\mu(e,e,e),\ B=\mu(e,e,f),\ C=\mu(e,f,f)$ and $D=\mu(f,f,f)$. From the exact sequence \eqref{SES: cup e isom} we see $H^{2}(N)\xrightarrow{\cup e}H^{4}(N)$ is an isomorphism. Its matrix representation with respect to the bases $(e,f)$ and $\left(e\spcheck,f\spcheck\right)$ is $M_{e}=\begin{pmatrix}
	A & B \\
	B & C
\end{pmatrix}$. Since this morphism is an isomorphism, we have $M_{e}\in\mathrm{GL}(2,\mathbb{Z})$. This justifies Condition 2 of Lemma \ref{Lemma: Base space 1}. We also obtain that $\left(e^{2},ef\right)$ is also a basis of $H^{4}(N)$.

Now we study $w_{2}(N)$ and $p_{1}(N)$. First we determine $w_{2}(N)$. Let $\xi$ be the $2-$plane bundle over $N$ associated to the circle bundle. Let $D_{e}\xrightarrow{\pi}N$ be the associated $D^{2}-$bundle. Let $M\xrightarrow{\iota}D_{e}$ be the boundary inclusion. Let $M\xrightarrow{p}N$ be the bundle projection as before. Then $p=\pi\circ\iota$. See the following commutative diagram.
\begin{align}\label{BundleDiagram}
	\xymatrix{
	 D_{e} \ar[r]^{\pi} & N \ar@{=}[d] \\
	 M \ar[r]^{p} \ar@{^(->}[u]^{\iota} & N 
	}
\end{align}
The stable tangent bundle of $M$ is $\tau_{M}\oplus\varepsilon^{1}=\iota^{*}\tau_{D_{e}}=\iota^{*}\pi^{*}\left(\tau_{N}\oplus\xi\right)=p^{*}\left(\tau_{N}\oplus\xi\right)$. Hence $$w_{2}\left(M\right)=p^{*}\left(w_{2}(N)+w_{2}(\xi)\right)=p^{*}w_{2}(N)+p^{*}\left(\overline{e}\right)=p^{*}w_{2}(N).$$ Here $\mathrm{mod}$ $2$ Gysin sequence implies $p^{*}\left(\overline{e}\right)=0$. While $M$ is spin, we have $w_{2}\left(M\right)=p^{*}w_{2}(N)=0$. Hence by $\mathrm{mod}$ $2$ Gysin sequence $w_{2}(N)=0$ or $\overline{e}$. Next we study $p_{1}(N)$. Set $p_{1}(N)=k\cdot e\spcheck+l\cdot f\spcheck$, $k,l\in\mathbb{Z}$. We begin with the relation 
\begin{eqnarray}\label{Formula: p1 and w2}
	p_{1}(N)\ \mathrm{mod}\ 2=w_{2}(N)^{2}.
\end{eqnarray}
This formula implies that 
\begin{compactenum}
	\item if $N$ is spin, $k\equiv l\equiv0\ \mathrm{mod}\ 2$;
	\item if $N$ is nonspin with $w_{2}(B)=\overline{e}$, $k\equiv A\ \mathrm{mod}\ 2$ and $l\equiv B\ \mathrm{mod}\ 2$.
\end{compactenum}
Now we relate the cohomology ring and characteristic classes of $N$. We begin with the Wu formula (\cite[\S 11, Theorem 11.14]{MilnorStasheff}). For our oriented $6-$manifold $N$ we have 
\begin{eqnarray}\label{Formula: Wu formula}
	\left<\mathrm{Sq^{2}}(x),[N]\right>=\left<w_{2}(N)\cup x,[N]\right>,\ \forall x\in H^{4}(N;\mathbb{Z}/2).
\end{eqnarray}
Hence by setting $x=\overline{e}^{2},\ \overline{e}\cdot\overline{f}$ and $\overline{f}^{2}$, we obtain that
\begin{compactenum}
	\item if $N$ is spin, $B+C\equiv0\ \mathrm{mod}\ 2$;
	\item if $N$ is nonspin with $w_{2}(N)=\overline{e}$, $A\equiv C\equiv0\ \mathrm{mod}\ 2$.
\end{compactenum}
Further restriction is given by \cite[Theorem 1]{Jupp1973} as follows. Let $\widehat{w}\in H^{2}(N)$ be an integral cohomology class that restricts to $w_{2}(N)\in H^{2}(N;\mathbb{Z}/2)$. Then we have the formula:
\begin{eqnarray}\label{Jupp relation}
	\mu\left(\widehat{w}+2x,\widehat{w}+2x,\widehat{w}+2x\right)\equiv\left<p_{1}(N)\cup\left(\widehat{w}+2x\right),[N]\right>\ \mathrm{mod}\ 48,\ \forall x\in H^{2}(N).
\end{eqnarray}
We take $\widehat{w}=0$ when $N$ is spin. In formula \eqref{Jupp relation} we set $x=Xe+Yf,\ X,Y\in\mathbb{Z}$, and we have
\begin{eqnarray}\label{Formula: Jupp relation, case 1, base spin}
	4X^{3}A+12X^{2}YB+12XY^{2}C+4Y^{3}D\equiv kX+lY\ \mathrm{mod}\ 24,\ \forall X,Y\in\mathbb{Z}.
\end{eqnarray}
Combining formulae \eqref{Formula: p1 and w2}, \eqref{Formula: Wu formula} and \eqref{Formula: Jupp relation, case 1, base spin}, we obtain:
\begin{eqnarray}\label{Formula: Jupp relation, case 1, base spin, expand}
	\begin{dcases}
		k \equiv 4A &\mod24,\\
		l \equiv 4D &\mod24,\\
		A \equiv 0 &\mod2,\\
		B \equiv C \equiv 1 &\mod2.
	\end{dcases}
\end{eqnarray}
When $N$ is nonspin with $w_{2}(N)=\overline{e}$, we take $\widehat{w}=e$. Set $x=Xe+Yf$, $\widehat{w}+2x=(2X+1)e+2Yf=\widetilde{X}e+2Yf$ with $\widetilde{X}=2X+1$, and formula \eqref{Jupp relation} becomes
\begin{eqnarray}\label{Formula: Jupp relation, case 1, base nonspin}
	\widetilde{X}^{3}A+6\widetilde{X}^{2}YB+12\widetilde{X}Y^{2}+8Y^{3}D\equiv k\widetilde{X}+2lY\ \mathrm{mod}\ 48,\ \forall\widetilde{X}\in\mathbb{Z}\ \text{odd}\ \forall Y\in\mathbb{Z}.
\end{eqnarray}
Combining formulae \eqref{Formula: p1 and w2}, \eqref{Formula: Wu formula} and \eqref{Formula: Jupp relation, case 1, base nonspin}, we obtain:
\begin{eqnarray}\label{Formula: Jupp relation, case 1, base nonspin, expand}
	\begin{dcases}
		k \equiv A & \mod48,\\
		l \equiv 3B+6C+4D & \mod24,\\
		A \equiv C \equiv 0 & \mod2,\\
		B \equiv 1 & \mod2.
	\end{dcases}
\end{eqnarray}
This justifies Condition 3 of Lemma \ref{Lemma: Base space 1} and hence we complete the proof of Lemma \ref{Lemma: Base space 1}.
\end{pf}
\begin{lem}\label{Lemma: Base space 1, if part}
	Let $A,\ B,\ C,\ D,\ u,\ v$ be $6$ integers such that $A$ is even and $B,C$ are odd (resp. $A,C$ are even and $B$ is odd). Up to orientation-preserving diffeomorphism there is a unique simply connected spin (resp. nonspin) $6-$manifold $N$ that satisfies Conditions 1, 2 and 3(a) (resp. Conditions 1, 2 and 3(b)) of Lemma \ref{Lemma: Base space 1} and is the orbit of a free circle action on a spin cohomology $S^{2}\times S^{5}$.
\end{lem}
\begin{pf}
It follows from \cite[Theorem 1]{Jupp1973} that up to orientation-preserving diffeomorphism there is a unique simply connected spin (resp. nonspin) $6-$manifold $N$ which satisfies Conditions 1, 2 and 3(a) (resp. Conditions 1, 2 and 3(b)) of Lemma \ref{Lemma: Base space 1}. Let $N_{e}$ be the total space of the circle bundle over $N$ whose Euler class is $e$. We will show that $N_{e}$ is a spin cohomology $S^{2}\times S^{5}$. 

Long exact sequence of homotopy groups associated to the circle bundle $N_{e}\xrightarrow{p}N$ implies $\pi_{1}\left(N_{e}\right)$ is a cyclic group. Hence $\pi_{1}\left(N_{e}\right)$ is abelian and by Hurewicz theorem $H_{1}\left(N_{e}\right)\cong\pi_{1}\left(N_{e}\right)$. Combining Gysin sequence and the fact that $e$ is primitive we have $H^{1}\left(N_{e}\right)=0$ and $H^{2}\left(N_{e}\right)\cong\mathbb{Z}$. Hence by universal coefficient theorem $H_{1}\left(N_{e}\right)=0$ and $N_{e}$ is simply connected. Since $H^{2}(N)\xrightarrow{\cup e}H^{4}(N)$ is an isomorphism, by Gysin sequence $H^{3}\left(N_{e}\right)=H^{4}\left(N_{e}\right)=0$. Now we see that $N_{e}$ has torsion-free cohomology groups with only nonzero Betti numbers $b_{0}=b_{2}=b_{5}=b_{7}=1$. Hence the total space $N_{e}$ is a cohomology $S^{2}\times S^{5}$. The formula $w_{2}\left(N_{e}\right)=p^{*}\left(w_{2}(N)+\overline{e}\right)$ still applies. By the assumption that $w_{2}(N)=0$ or $\overline{e}$ and the mod $2$ Gysin sequence we have $w_{2}\left(N_{e}\right)$. Hence $N_{e}$ is spin. 
\end{pf}

\subsection{Invariants of $N_{e}$}\label{Subsec: s inv 1}
In this section we compute the homeomorphism and diffeomorphism invariants of $N_{e}$, where $(N,e)$ is given as in Lemma \ref{Lemma: Base space 1}. The homeomorphism and diffeomorphism invariants are the $s-$invariants developed by Kreck and Stolz (\cite{KS88}, \cite{Kreck1991SomeNH} and \cite{Kreck1991SomeNHCorrection}). We first review the definition of $s-$invariants and then compute the $s-$invariants of the manifold $N_{e}$.

In \cite{KS88}, \cite{Kreck1991SomeNH} and \cite{Kreck1991SomeNHCorrection} the authors considered the homeomorphism and diffeomorphism classification of the $7-$manifolds satisfying the following condition (\cite[Condition 2.1]{Kreck1991SomeNH}):
\begin{quotation}
	The $7-$manifold $M$ is simply connected. Its second cohomology group is infinite cyclic with a generator $z$, and its forth cohomology group is finite cyclic with the generator $z^{2}$.
\end{quotation}
For convenience a spin manifold satisfying the condition above is called a \emph{manifold of type I}. In particular any spin cohomology $S^{2}\times S^{5}$ is a type I manifold. By assumption $N_{e}$ is a spin cohomology $S^{2}\times S^{5}$ and thus a type I manifold with $z=p^{*}f$ a generator of the second cohomology.

The homeomorphism and diffeomorphism invariants of the type I manifold M are defined as follows. Let $W$ be a compact $8-$manifold whose boundary is $M$. We also say $W$ is a coboundary of $M$. Let $\widehat{z}\in H^{2}(W)$ restrict to $z$. When $W$ is spin the pair $\left(W,\widehat{z}\right)$ is called a spin coboundary of $(M,z)$. By \cite[Lemma 6.1]{Kreck1991SomeNH} such a spin coboundary of $(M,z)$ always exists. Sometimes we cannot easily construct such a spin coboundary, but we can find a nonspin coboundary $W$ and its cohomology classes $\widehat{z},\ c\in H^{2}(W)$ such that $w_{2}(W)=\overline{c}$ and $\widehat{z}|_{M}=z,\ c|_{M}=0$. In this case we say the triple $\left(W,\widehat{z},c\right)$ is a nonspin coboundary of $(M,z)$. Then we have the following characteristic numbers (\cite[Formula 2.7]{Kreck1991SomeNH}):
\begin{eqnarray}\label{Formula: S invariant}
	\left\{\begin{aligned}
		S_{1}\left(W,\widehat{z},c\right)& =\left<e^{\frac{c}{2}}\widehat{A}(W),[W,M]\right>\in\mathbb{Q},\\
		S_{2}\left(W,\widehat{z},c\right)& =\left<\left(e^{\widehat{z}}-1\right)e^{\frac{c}{2}}\widehat{A}(W),[W,M]\right>\in\mathbb{Q},\\
		S_{3}\left(W,\widehat{z},c\right)& =\left<\left(e^{2\widehat{z}}-1\right)e^{\frac{c}{2}}\widehat{A}(W),[W,M]\right>\in\mathbb{Q},
	\end{aligned}\right.
\end{eqnarray}
where we take $c=0$ if $\left(W,\widehat{z}\right)$ is a spin coboundary of $(M,z)$. The terms on the right hand sides are understood as follows. The expression $\widehat{A}(W)$ is the $\widehat{A}-$polynomial of $W$. The cohomology classes $\left(e^{\widehat{z}}-1\right)e^{\frac{c}{2}}\widehat{A}(W)$ and $\left(e^{2\widehat{z}}-1\right)e^{\frac{c}{2}}\widehat{A}(W)$ are rational combinations of $p_{1}^{2}$, $\widehat{z}^{i}c^{j}p_{1}(i+j=2)$ and $\widehat{z}^{i}c^{j}(i+j=4)$. Here $p_{1}=p_{1}(W)$. By assumption the classes $p_{1}$ and $\widehat{z}^{i}c^{j}(i+j=2)$ restrict to $0\in H^{4}(M;\mathbb{Q})$, hence they can be viewed as classes in $H^{4}(W,M;\mathbb{Q})$. Then the classes $p_{1}^{2}$, $\widehat{z}^{i}c^{j}p_{1}(i+j=2)$ and $\widehat{z}^{i}c^{j}(i+j=4)$ can be viewed as classes in $H^{8}(W,M;\mathbb{Q})$ and evaluate on $[W,M]$. The cohomology class $e^{\frac{c}{2}}\widehat{A}(W)$ is a rational linear combination of $p_{1}^{2}$, $\widehat{z}^{i}c^{j}p_{1}(i+j=2)$, $\widehat{z}^{i}c^{j}(i+j=4)$ and the $L-$polynomial $L(W)=\frac{1}{45}\left(7p_{2}-p_{1}^{2}\right)$. The term $\left<L(W),[W,M]\right>$ is interpreted as the signature $\sigma(W,M)$ of the pair $(W,M)$, namely the signature of the intersection form on $H^{4}(W,M;\mathbb{Q})$. 

We set $s_{i}(M)=S_{i}\left(W,\widehat{z},c\right)\ \mathrm{mod}\ 1\in\mathbb{Q}/\mathbb{Z}$, $i=1,2,3$. These $s_{i}-$invariants do not depend on the choices of the generators of $H^{2}(M)$ or the coboundaries $\left(W,\widehat{z},c\right)$. Hence the $s_{i}-$invariants are invariants of the type I manifold $M$. It follows from \cite[Theorem 1]{Kreck1991SomeNHCorrection} that two manifolds $M$, $M'$ of type I are diffeomorphic (resp. homeomorphic) if and only if their forth cohomology groups have the same order and they share the same invariants $s_{1}$, $s_{2}$ and $s_{3}$ (resp. $28s_{1}$, $s_{2}$ and $s_{3}$).

With a slight abuse of notation we write $p_{1}^{2}$ for $\left<p_{1}^{2},[W,M]\right>$, and the other monomials are similar. We have the following explicit expansion of formula \eqref{Formula: S invariant}.
\begin{eqnarray}\label{Formula: s invariant, expand}
	\left\{\begin{aligned}
		S_{1}\left(W,\widehat{z},c\right)& = -\frac{1}{224}\sigma(W,M)+\frac{1}{896}p_{1}^{2}-\frac{1}{192}c^{2}p_{1}+\frac{1}{384}c^{4},\\
		S_{2}\left(W,\widehat{z},c\right)& =-\frac{1}{48}\left(\widehat{z}^{2}+\widehat{z}c\right)p_{1}+\frac{1}{48}\left(2\widehat{z}^{4}+4\widehat{z}^{3}c+3\widehat{z}^{2}c^{2}+\widehat{z}c^{3}\right),\\
		S_{3}\left(W,\widehat{z},c\right)& =-\frac{1}{24}\left(2\widehat{z}^{2}+zc\right)p_{1}+\frac{1}{24}\left(16\widehat{z}^{4}+16\widehat{z}^{3}c+6\widehat{z}^{2}c^{2}+\widehat{z}c^{3}\right).
	\end{aligned}\right.
\end{eqnarray}
\begin{eg}
Let $(N,e)$ be given as in Lemma \ref{Lemma: Base space 1}. When $N$ is spin and the matrix $M_{e}$ has determinant $-1$, the $s-$invariants of the type I manifold $N_{e}$ are given as follows:
\begin{eqnarray}\label{eq:s^I inv of N_e, N spin and det M_e=-1 }
	\left\{\begin{aligned}
		s_{1}\left(N_{e}\right)&=-\frac{9}{14}\left(Cu^{2}-2Buv+Av^{2}\right)+\frac{2-3B(B-D)}{14}u+\frac{3A(B-D)}{14}v+\frac{A}{224}-\frac{A}{56}(B-D)^{2}\ \mathrm{mod}\ 1,\\
		28s_{1}\left(N_{e}\right)&=\frac{A}{8}\ \mathrm{mod}\ 1;\\
		s_{2}\left(N_{e}\right)&=\frac{1+D}{2}u+\frac{2AC^{2}-2ABD+AD^{2}}{24}\ \mathrm{mod}\ 1;\\
		s_{3}\left(N_{e}\right)&=\frac{1}{3}\left(AC^{2}-ABD-AD^{2}\right)+\frac{D+1}{2}\ \mathrm{mod}\ 1.
	\end{aligned}
	\right.
\end{eqnarray}
And when $N$ is nonspin with $w_{2}(N)=\overline{e}$, the $s-$invariants of the type I manifold $N_{e}$ are given as follows:
\begin{eqnarray}\label{eq:s^I inv of N_e, N nonspin}
	\left\{\begin{aligned}
		s_{1}\left(N_{e}\right)&=-\frac{9}{14}\left(4Cu^{2}-4Buv+Av^{2}\right)+\frac{3}{14}\left(B^{2}+3BC+2BD\right)u-\frac{3}{28}(AB+3AC+2AD)v\\
		&\ -\frac{1}{224}\left(A^{2}C+6ABC+4ABD+9AC^{2}+4AD^{2}+12ACD\right)\ \mathrm{mod}\ 1,\\
		28s_{1}\left(N_{e}\right)&=0\ \mathrm{mod}\ 1;\\
		s_{2}\left(N_{e}\right)&=-\frac{1}{24}\left(B^{2}C+3BC^{2}-ABD-3ACD-AD^{2}+C+C^{3}\right)\ \mathrm{mod}\ 1;\\
		s_{3}\left(N_{e}\right)&=-\frac{1}{6}\left(B^{2}C-ABD+2AD^{2}+C+C^{3}\right)\ \mathrm{mod}\ 1.
	\end{aligned}
	\right.
\end{eqnarray}
When $N$ is spin we omit the discussion of the case $\mathrm{det}\ M_{e}=1$, since it turns out that the case $\mathrm{det}\ M_{e}=-1$ is sufficient. It is also useful to notice that when $N$ is nonpsin with $w_{2}(N)=\overline{e}$ we have $s_{3}\left(N_{e}\right)=4s_{2}\left(N_{e}\right)$.

The $s-$invariants of the manifold $N_{e}$ are computed as follows. First we need to find a coboundary of $N_{e}$. A natural choice is the associated disc bundle $D_{e}$. From the Gysin sequence and Thom isomorphism we have the following commutative diagram.
\begin{eqnarray}\label{diagram: Thom isom and Gysin seq}
	\xymatrix{
		0\ar[r] & H^{2}\left(D_{e},N_{e}\right)\ar^{j^{*}}[r] & H^{2}\left(D_{e}\right)\ar[r]^{\iota^{*}} & H^{2}\left(N_{e}\right)\ar@{=}[d]\ar[r] & 0\\
		0\ar[r] & H^{0}(N)\ar[r]^{\cup e}\ar[u]^{\cup u(\xi)}_{\cong} & H^{2}(N)\ar[r]^{p^{*}}\ar[u]^{\pi^{*}}_{\cong} & H^{2}\left(N_{e}\right)\ar[r] & 0
	}
\end{eqnarray}
Here $j:\left(D_{e},\varnothing\right)\to\left(D_{e},N_{e}\right)$ is the inclusion and the Thom class of $\xi$ is $u(\xi)\in H^{2}\left(D_{e},N_{e}\right)$. In diagram \eqref{diagram: Thom isom and Gysin seq} the rows are exact and the vertical morphisms are isomorphisms or identity. From the diagram we obtain that the boundary inclusion $N_{e}\xrightarrow{\iota}D_{e}$ induces an epimorphism on second cohomology groups, such that $H^{2}\left(D_{e}\right)=\mathbb{Z}\left\{\pi^{*}e,\pi^{*}f\right\}\cong\mathbb{Z}^{2}$, $\iota^{*}\left(\pi^{*}e\right)=0$ and $\iota^{*}\left(\pi^{*}f\right)=p^{*}f$ generates $H^{2}\left(N_{e}\right)\cong\mathbb{Z}$. From the formula $w_{2}\left(D_{e}\right)=\pi^{*}\left(w_{2}\left(N\right)+w_{2}(\xi)\right)=\pi^{*}\left(w_{2}\left(N\right)+\overline{e}\right)$ we see that $D_{e}$ is spin if $N$ is not and nonspin if $N$ is. If we keep using natural coboundary $W=D_{e}$ of $N_{e}$, we should treat these two cases separately. When $D_{e}$ is spin we use the spin coboundary $\left(D_{e},\pi^{*}f\right)$, and when $D_{e}$ is nonspin we use the nonspin coboundary $\left(D_{e},\pi^{*}f,\pi^{*}e\right)$.

Next we compute the terms on the right hand sides of formula \eqref{Formula: s invariant, expand}. We begin with the signature term and claim that $\sigma\left(D_{e},N_{e}\right)=\sigma\left(M_{e}\right)$. This follows from the following commutative diagram
\begin{eqnarray}\label{Diagram: bilinear form isom}
	\xymatrix{
		H^{4}\left(D_{e},N_{e}\right)\otimes H^{4}\left(D_{e},N_{e}\right)\ar[r]^(0.525){\mathrm{id}\otimes j^{*}}_{\cong} & 	H^{4}\left(D_{e},N_{e}\right)\otimes H^{4}\left(D_{e}\right)\ar[r]^(0.6){\cup} & H^{8}\left(D_{e},N_{e}\right) \\
		H^{2}(N)\otimes H^{2}(N)\ar[u]^{\cup u(\xi)\otimes\cup u(\xi)}_{\cong}\ar[r]^{\mathrm{id}\otimes\cup e}_{\cong} & H^{2}(N)\otimes H^{4}(N)\ar[u]^{\cup u(\xi)\otimes\pi^{*}}_{\cong}\ar[r]^(0.6){\cup} & H^{6}(N)\ar[u]^{\cup u(\xi)}_{\cong}
	}
\end{eqnarray}
Here the first two vertical morphisms are isomorphisms, since they are tensors of isomorphisms between finitely generated free abelian groups. Hence the bilinear forms $\left(H^{4}\left(D_{e},N_{e}\right),\cdot\cup\cdot\right)$ and $\left(H^{2}(N),\left(\cdot\cup\cdot\right)\cup e\right)$ are isomorphic. Therefore $\sigma\left(D_{e},N_{e}\right)=\sigma\left(M_{e}\right)$. We denote this common value by $\sigma$. Next we compute the products of $p_{1}$, $\widehat{z}$ and $c$. From diagram \eqref{Diagram: bilinear form isom} we observe the following helpful fact. Given $\widetilde{u}$ and $\widetilde{v}\in H^{4}\left(D_{e}\right)$, there are unique classes $u$ and $v\in H^{2}(N)$ such that $\widetilde{u}=j^{*}\left(u\cup u(\xi)\right)=\pi^{*}\left(e\cup u\right)$ and $\widetilde{v}=\pi^{*}\left(e\cup v\right)$. Then $\widetilde{u}\cdot\widetilde{v}=\mu(e,u,v)$. Note $c^{2}=\pi^{*}\left(e\cup e\right)$, $\widehat{z}c=\pi^{*}\left(e\cup f\right)$, and $f^{2}=\begin{pmatrix}
	e\spcheck & f\spcheck
\end{pmatrix}\begin{pmatrix}
	C \\ D
\end{pmatrix}=e\cup\begin{pmatrix}
	e & f
\end{pmatrix}M_{e}^{-1}\begin{pmatrix}
	C \\ D
\end{pmatrix}$. Hence $\widehat{z}^{2}=\pi^{*}\left(f^{2}\right)=\pi^{*}\left(e\cup\begin{pmatrix}
	e & f
\end{pmatrix}M_{e}^{-1}\begin{pmatrix}
	C \\ D
\end{pmatrix}\right).$ From $\tau_{D_{e}}\cong\pi^{*}\left(\tau_{N}\oplus\xi\right)$ we have $$p_{1}\left(D_{e}\right)=\pi^{*}\left(p_{1}(N)+e^{2}\right)=\pi^{*}\left(e\cup\left(\begin{pmatrix}
e & f
\end{pmatrix}M_{e}^{-1}\begin{pmatrix}
k \\ l
\end{pmatrix}+e\right)\right).$$
The characteristic numbers are given as follows.
\begin{align*}
		p_{1}^{2}=p_{1}\cdot p_{1} & = \mu\left(e,\begin{pmatrix}
		e & f
	\end{pmatrix}M_{e}^{-1}\begin{pmatrix}
		k \\ l
	\end{pmatrix}+e,\begin{pmatrix}
		e & f
	\end{pmatrix}M_{e}^{-1}\begin{pmatrix}
		k \\ l
	\end{pmatrix}+e\right)=A+2k+\begin{pmatrix}
	k & l
	\end{pmatrix}M_{e}^{-1}\begin{pmatrix}
	k \\ l
	\end{pmatrix};
\end{align*}
\begin{align*}
	c^{2}\cdot p_{1} & =\mu\left(e,e,\begin{pmatrix}
		e & f
	\end{pmatrix}M_{e}^{-1}\begin{pmatrix}
		k \\ l
	\end{pmatrix}+e\right)
%	=\begin{pmatrix}
%		1 & 0
%	\end{pmatrix}M_{e}M_{e}^{-1}\begin{pmatrix}
%		k \\ l
%	\end{pmatrix}+A
	=k+A,\\
	zc\cdot p_{1} & =\mu\left(e,f,\begin{pmatrix}
		e & f
	\end{pmatrix}M_{e}^{-1}\begin{pmatrix}
		k \\ l
	\end{pmatrix}+e\right)
%	=\begin{pmatrix}
%		0 & 1
%	\end{pmatrix}M_{e}M_{e}^{-1}\begin{pmatrix}
%		k \\ l
%	\end{pmatrix}+B
	=l+B,\\
	z^{2}\cdot p_{1} & = \mu\left(e,\begin{pmatrix}
		e & f
	\end{pmatrix}M_{e}^{-1}\begin{pmatrix}
		C \\ D
	\end{pmatrix},\begin{pmatrix}
		e & f
	\end{pmatrix}M_{e}^{-1}\begin{pmatrix}
		k \\ l
	\end{pmatrix}+e\right)%\\
%	& = \begin{pmatrix}
%		1 & 0
%	\end{pmatrix}M_{e}M_{e}^{-1}\begin{pmatrix}
%		C \\ D
%	\end{pmatrix}+\begin{pmatrix}
%		C & D
%	\end{pmatrix}\left(M_{e}^{-1}\right)^{T}M_{e}M_{e}^{-1}\begin{pmatrix}
%		k \\ l
%	\end{pmatrix}\\ &
	 = C+\begin{pmatrix}
		C & D
	\end{pmatrix}M_{e}^{-1}\begin{pmatrix}
		k \\ l
	\end{pmatrix};\\
%\end{align*}
%\begin{align*}
	c^{2}\cdot c^{2} & =A,\\
	\widehat{z}c\cdot c^{2} & =B,\\
	\widehat{z}c\cdot \widehat{z}c & =C,\\
	\widehat{z}^{2}\cdot c^{2} & =\mu\left(e,\begin{pmatrix}
		e & f
	\end{pmatrix}M_{e}^{-1}\begin{pmatrix}
		C \\ D
	\end{pmatrix},e\right)=C,\\
	\widehat{z}^{2}\cdot \widehat{z}c & =\mu\left(e,\begin{pmatrix}
		e & f
	\end{pmatrix}M_{e}^{-1}\begin{pmatrix}
		C \\ D
	\end{pmatrix},f\right)=D,\\
	\widehat{z}^{2}\cdot \widehat{z}^{2} & = \mu\left(e,\begin{pmatrix}
		e & f
	\end{pmatrix}M_{e}^{-1}\begin{pmatrix}
		C \\ D
	\end{pmatrix},\begin{pmatrix}
		e & f
	\end{pmatrix}M_{e}^{-1}\begin{pmatrix}
		C \\ D
	\end{pmatrix}\right) = \begin{pmatrix}
		C & D
	\end{pmatrix}M_{e}^{-1}\begin{pmatrix}
		C \\ D
	\end{pmatrix}.
\end{align*}
Here we have two ways to compute $\widehat{z}^{2}c^{2}$, and they lead the same result.

Now we express the $s-$invariants of $N_{e}$ explicitly. When $N$ is spin we have $w_{2}\left(D_{e}\right)=\pi^{*}\left(\overline{e}\right)=\left(\pi^{*}e\right)\ \mathrm{mod}\ 2$, $\left.\left(\pi^{*}e\right)\right|_{N_{e}}=0$ and $\left.\left(\pi^{*}f\right)\right|_{N_{e}}=p^{*}f$. Since $\mathrm{det}M_{e}=-1$, we have $M_{e}^{-1}=\begin{pmatrix}
-C & B\\ B & -A
\end{pmatrix}$ and $\sigma=0$. The $S-$invariants are
\begin{align*}
	S_{1}\left(D_{e},z,c\right) & =\frac{1}{896}\left(A+2k+\begin{pmatrix}
		k & l
	\end{pmatrix}M_{e}^{-1}\begin{pmatrix}
		k \\ l
	\end{pmatrix}\right)-\frac{1}{192}(k+A)+\frac{1}{384}A\\
%	& = \frac{1}{896}\left(A+2k-Ck^{2}+2Bkl-Al^{2}\right)-\frac{1}{192}(k+A)+\frac{1}{384}A\\
%	& = -\frac{1}{672}A-\frac{1}{336}k-\frac{1}{896}\left(Ck^{2}-2Bkl+Al^{2}\right)\\
%	& = -\frac{1}{672}A-\frac{1}{336}(24u+4A)-\frac{1}{896}\left(C(24u+4A)^{2}-2B(24u+4A)(24v+4D)+A(24v+4D)^{2}\right)\\
	& = -\frac{9}{14}\left(Cu^{2}-2Buv+Av^{2}\right)+\frac{2-3B(B-D)}{14}u+\frac{3A(B-D)}{14}v+\frac{A}{224}-\frac{A}{56}(B-D)^{2},\\
%\end{align*}
%\begin{align*}
	S_{2}\left(D_{e},z,c\right) & = -\frac{1}{48}\left(C+\begin{pmatrix}
		C & D
	\end{pmatrix}M_{e}^{-1}\begin{pmatrix}
		k \\ l
	\end{pmatrix}+(l+B)\right)+\frac{1}{48}\left(2\begin{pmatrix}
	C & D
	\end{pmatrix}M_{e}^{-1}\begin{pmatrix}
	C \\ D
	\end{pmatrix}+4D+3C+B\right)\\
%	& = \frac{1}{48}\left(2C+C^{2}k-B(Cl+Dk)+ADl-l-2C^{3}+4BCD-2AD^{2}+4D\right)\\
%	& = \frac{1}{48}\left(\left(C^{2}-BD\right)k+(AD-BC-1)l+2C-2C^{3}+4BCD-2AD^{2}+4D\right)\\
%	& = \frac{1}{48}\left(\left(C^{2}-BD\right)(24u+4A)+(AD-BC-1)(24v+4D)+2C-2C^{3}+4BCD-2AD^{2}+4D\right)\\
	& = \frac{C^{2}-BD}{2}u+\frac{AD-BC-1}{2}v+\frac{2AC^{2}-2ABD+AD^{2}+C-C^{3}}{24},\\
%\end{align*}
%\begin{align*}
	S_{3}\left(D_{e},z,c\right) & = -\frac{1}{24}\left(2\left(C+\begin{pmatrix}
		C & D
	\end{pmatrix}M_{e}^{-1}\begin{pmatrix}
		k \\ l
	\end{pmatrix}\right)+(l+B)\right)+\frac{1}{24}\left(16\begin{pmatrix}
	C & D
	\end{pmatrix}M_{e}^{-1}\begin{pmatrix}
	C \\ D
	\end{pmatrix}+16D+6C+B\right)\\
%	& =\frac{1}{24}\left(4C+2C^{2}k-2B(Cl+Dk)+2ADl-l-16C^{3}+32BCD-16AD^{2}+16D\right)\\
%	& = \frac{1}{24}\left(2\left(C^{2}-BD\right)k+(2AD-2BC-1)l+4C-16C^{3}+32BCD-16AD^{2}+16D\right)\\
%	& = \frac{1}{24}\left(2\left(C^{2}-BD\right)(24u+4A)+(2AD-2BC-1)(24v+4D)+4C-16C^{3}+32BCD-16AD^{2}+16D\right)\\
	& = 2\left(C^{2}-BD\right)u+(2AD-2BC-1)v\\
	&\quad +\frac{1}{24}\left(8AC^{2}-8ABD-8AD^{2}+24BCD+12D+4C-16C^{3}\right).
\end{align*}
By Lemma \ref{Lemma: Base space 1}, Condition 3(a), the parameter $A$ is even and $B$ and $C$ are odd. Hence modulo $\mathbb{Z}$ we obtain the formula \eqref{eq:s^I inv of N_e, N spin and det M_e=-1 }. When $N$ is nonspin with $w_{2}(N)=\overline{e}$ we have $w_{2}\left(D_{e}\right)=0$. Take $c=0$ in \eqref{Formula: S invariant}. Recall $\sigma=0$ and $\det M_{e}=AC-B^{2}=-1<0$. The $S-$invariants are
\begin{align*}
	S_{1}\left(D_{e},z\right)% & = \left<\widehat{A}\left(D_{e}\right),\left[D_{e},B_{e}\right]\right>= \frac{1}{896}p_{1}^{2}\\
%	& =\frac{1}{896}\left(A+2k-Ck^{2}+2Bkl-Al^{2}\right)\\
%	& = \frac{1}{896}\left(A+2(48u+A)-C(48u+A)^{2}+2B(48u+A)(24v+3B+6C+4D)\right.\\
%	& \ \left.-A(24v+3B+6C+4D)^{2}\right)\\
	& =-\frac{9}{14}\left(4Cu^{2}-4Buv+Av^{2}\right)+\frac{3}{14}\left(B^{2}+3BC+2BD\right)u-\frac{3}{28}(AB+3AC+2AD)v\\
	&\ -\frac{1}{224}\left(A^{2}C+6ABC+4ABD+9AC^{2}+4AD^{2}+12ACD\right),\\
%\end{align*}
%\begin{align*}
	S_{2}\left(D_{e},z\right) %& = \left<\mathrm{ch}\left(\lambda(z)-1\right)\widehat{A}\left(D_{e}\right),\left[D_{e},B_{e}\right]\right>=-\frac{1}{48}z^{2}p_{1}+\frac{1}{24}z^{4}\\
%	& = -\frac{1}{48}\left(C-C^{2}k+B(Cl+Dk)-ADl\right)-\frac{1}{24}\left(C^{3}-2BCD+AD^{2}\right)\\
%	& = \frac{1}{48}\left(\left(C^{2}-BD\right)k-(BC-AD)l\right)-\frac{1}{48}C-\frac{1}{24}\left(C^{3}-2BCD+AD^{2}\right)\\
%	& =\frac{1}{48}\left(\left(C^{2}-BD\right)(48u+A)-(BC-AD)(24v+3B+6C+4D)\right)-\frac{1}{48}C-\frac{1}{24}\left(C^{3}-2BCD+AD^{2}\right)\\
	& = \left(C^{2}-BD\right)u+\frac{AD-BC}{2}v-\frac{1}{24}\left(B^{2}C+3BC^{2}-ABD-3ACD-AD^{2}+C+C^{3}\right),\\
%\end{align*}
%\begin{align*}
	S_{3}\left(D_{e},z\right)% & = \left<\mathrm{ch}\left(\lambda^{2}(z)-1\right)\widehat{A}\left(D_{e}\right),\left[D_{e},B_{e}\right]\right>= -\frac{1}{12}z^{2}p_{1}+\frac{2}{3}z^{4}\\
%	& =-\frac{1}{12}\left(C-C^{2}k+B(Cl+Dk)-ADl\right)-\frac{2}{3}\left(C^{3}-2BCD+AD^{2}\right)\\
%	& = \frac{1}{12}\left(\left(C^{2}-BD\right)k+(AD-BC)l\right)-\frac{1}{12}C-\frac{2}{3}\left(C^{3}-2BCD+AD^{2}\right)\\
%	& = \frac{1}{12}\left(\left(C^{2}-BD\right)(48u+A)+(AD-BC)(24v+3B+6C+4D)\right)-\frac{1}{12}C-\frac{2}{3}\left(C^{3}-2BCD+AD^{2}\right)\\
	& = 4\left(C^{2}-BD\right)u+2(AD-BC)v-\frac{1}{6}\left(B^{2}C+3BC^{2}-ABD-3ACD+2AD^{2}+C+C^{3}\right).
\end{align*}
By Lemma \ref{Lemma: Base space 1}, Condition 3(b), the parameters $A$ and $C$ are even and $B$ is odd. Hence modulo $\mathbb{Z}$ we obtain formula \eqref{eq:s^I inv of N_e, N nonspin}.
\end{eg}

\subsection{Comparing $s-$invariants}\label{Subsec: compare s inv 1}
In this section we first compute the $s-$invariants of $S^{2}\times S^{5}\#\Sigma_{r}$. Then we compare the $s-$invariants of $S^{2}\times S^{5}\#\Sigma_{r}$ and the manifold $N_{e}$, deducing for which homotopy $7-$sphere $\Sigma$ the manifold $S^{2}\times S^{5}\#\Sigma$ can be written as the form $N_{e}$ and thus admits a free circle action.

First we compute the $s-$invariants of $S^{2}\times S^{5}\#\Sigma_{r}$. It follows from \cite[Section 3]{EKinv} that the $\mu-$invariant is defined for $S^{2}\times S^{5}$ and $\Sigma_{r}$, hence the $\mu-$invariant is also defined for their connected sum and $$\mu\left(S^{2}\times S^{5}\#\Sigma_{r}\right)=\mu\left(S^{2}\times S^{5}\right)+\mu\left(\Sigma_{r}\right).$$ We can use the coboundary $S^{2}\times D^{6}$ of $S^{2}\times S^{5}$ to compute that $\mu\left(S^{2}\times S^{5}\right)=0\ \mathrm{mod}\ 1$. Hence $$\mu\left(S^{2}\times S^{5}\#\Sigma_{r}\right)=\mu\left(\Sigma_{r}\right)=\frac{r}{28}\ \mathrm{mod}\ 1.$$ According to \cite[Section 3]{KS88}, if the invariants $s_{1}$ and $\mu$ are both defined for the type I manifold $M$, then $s_{1}(M)=\mu(M)$. Hence $s_{1}\left(S^{2}\times S^{5}\#\Sigma\right)=\mu\left(S^{2}\times S^{5}\#\Sigma_{r}\right)=\frac{r}{28}\ \mathrm{mod}\ 1$ and $28s_{1}\left(S^{2}\times S^{5}\#\Sigma_{r}\right)=0\ \mathrm{mod}\ 1$. By \cite[Section 4]{KervaireMilnor63} each homotopy $7-$sphere bounds a parallelizable compact $8-$manifold, saying $W_{r}$ is a parallelizable compact coboundary of $\Sigma_{r}$. The boundary connected sum $\left(S^{2}\times D^{6}\right)\natural W_{r}$ is a coboundary of $S^{2}\times S^{5}\#\Sigma_{r}$ and can be applied to compute the invariants $s_{2}$ and $s_{3}$ of $S^{2}\times S^{5}\#\Sigma_{r}$. The result is $s_{2}\left(S^{2}\times S^{5}\#\Sigma_{r}\right)=s_{3}\left(S^{2}\times S^{5}\#\Sigma_{r}\right)=0\ \mathrm{mod}\ 1$.
\begin{prop}\label{Propotision: homeom S2 times S5}
	\begin{compactenum}
		\item The manifold $S^{2}\times S^{5}\#\Sigma$ admits a free circle action with spin orbit for \emph{all} homotopy $7-$sphere $\Sigma$.
		\item The manifold $S^{2}\times S^{5}\#\Sigma_{r}$ admits a free circle action with nonspin orbit if and only if $r$ is even. 
	\end{compactenum}
\end{prop}
In the proof we will use the following lemma from elementary number theory. Its proof is postponed to the end of this section.
\begin{lem}\label{Lemma from number theory}
	Fix a prime $p$ and integers $a$, $b$, $c$ with $\gcd(ab,p)=1$. Let $N(a,b,c;p)$ be the number of solutions to the equation $ax^{2}+by^{2}\equiv c\ \mathrm{mod}\ p$.
	Then $$N(a,b,c;p)=\begin{dcases}
		p+\left(\frac{-ab}{p}\right)(p-1), & \text{if } c\equiv0\ \mathrm{mod}\ p,\\
		p-\left(\frac{-ab}{p}\right), & \text{otherwise},
	\end{dcases}$$
	where $\left(\frac{\cdot}{p}\right)$ is the Legendre symbol. In particular, any element in $\mathbb{Z}/p$ is represented by the quadratic form $ax^{2}+by^{2}$.
\end{lem}
\begin{pf}[of Proposition \ref{Propotision: homeom S2 times S5}]
	Let $(N,e)$ be given as in Lemma \ref{Lemma: Base space 1}. The strategy is as follows. First we solve the equation
\begin{eqnarray}\label{Eq: which N_e is homeom S2 times S5}
28s_{1}\left(N_{e}\right)=s_{2}\left(N_{e}\right)=s_{3}\left(N_{e}\right)=0\ \mathrm{mod}\ 1.
\end{eqnarray}
The manifold $N_{e}$ is homeomorphic to $S^{2}\times S^{5}$ if and only if it is a solution to equation \eqref{Eq: which N_e is homeom S2 times S5}. Then for such a manifold $N_{e}$ that is homeomorphic to $S^{2}\times S^{5}$, its $s_{1}-$invariant must be of the form $s_{1}\left(N_{e}\right)=\frac{r}{28}\ \mathrm{mod}\ 1$ for some $r\in\mathbb{Z}/28$. Set
\begin{align*}
	\mathcal{R}&=\left\{r\in\mathbb{Z}/28\middle|N_{e}\text{ is homeomorphic to }S^{2}\times S^{5}\text{ and }s_{1}\left(N_{e}\right)=\frac{r}{28}\ \mathrm{mod}\ 1\right\},\\
	\mathcal{R}^{+}&=\left\{r\in\mathbb{Z}/28\middle|N\text{ is spin, }N_{e}\text{ is homeomorphic to }S^{2}\times S^{5}\text{ and }s_{1}\left(N_{e}\right)=\frac{r}{28}\ \mathrm{mod}\ 1\right\},\\
	\mathcal{R}^{-}&=\left\{r\in\mathbb{Z}/28\middle|N\text{ is nonspin, }N_{e}\text{ is homeomorphic to }S^{2}\times S^{5}\text{ and }s_{1}\left(N_{e}\right)=\frac{r}{28}\ \mathrm{mod}\ 1\right\},
\end{align*}
and it is clear that $\mathbb{Z}/28\supset\mathcal{R}\supset\mathcal{R}^{\pm}$. Then the manifold $S^{2}\times S^{5}\#\Sigma_{r}$ admits a free circle action (resp. with spin orbit, with nonspin orbit) if and only if $r\in\mathcal{R}$ (resp. $\mathcal{R}^{+}$, $\mathcal{R}^{-}$). We will show that $\mathcal{R}^{-}=2\mathbb{Z}/28$ and $\mathcal{R}^{+}=\mathbb{Z}/28$. In particular $\mathcal{R}=\mathbb{Z}/28$ and the manifold $S^{2}\times S^{5}\#\Sigma$ admits a free circle action for \emph{any} homotopy $7-$sphere $\Sigma$.

We begin with the case $N$ is nonspin, then the case $N$ is spin with $\det M_{e}=-1$. %It turns out that these two cases are sufficient, and we omit the case $B$ is spin with $\det M_{e}=1$.
When $N$ is nonspin, we already have $28s_{1}\left(N_{e}\right)=0\ \mathrm{mod}\ 1$. As $s_{3}\left(N_{e}\right)=4s_{2}\left(N_{e}\right)$, it suffices to solve $s_{2}\left(N_{e}\right)=0\ \mathrm{mod}\ 1$. Recall that $A$, $C$ are even and $B$ is odd. Set $A=2A_{1}$, $C=2C_{1}$, and we solve the following equation set
\begin{equation*}
	\begin{dcases}
		4A_{1}C_{1}-B^{2}&=-1,\\
		B^{2}C_{1}+6BC_{1}^{2}-A_{1}BD-6A_{1}C_{1}D+C_{1}+4C_{1}^{3}&\equiv 0\ \mathrm{mod}\ 12.
	\end{dcases}
\end{equation*}
By Chinese remainder theorem, the second equation is equivalent to the congruence equation set\\
\begin{equation*}
	\begin{dcases}
		B^{2}C_{1}+6BC_{1}^{2}-A_{1}BD-6A_{1}C_{1}D+C_{1}+4C_{1}^{3}\equiv 0&\mathrm{mod}\ 3,\\
		B^{2}C_{1}+6BC_{1}^{2}-A_{1}BD-6A_{1}C_{1}D+C_{1}+4C_{1}^{3}\equiv 0&\mathrm{mod}\ 4.
	\end{dcases}
\end{equation*}
thus $A_{1}$, $B$, $C_{1}$, $D$ are subject to the following equation set
\begin{equation}%\label{B nonspin ABCD}
	\begin{dcases}
		\left(B^{2}-1\right)C_{1}\equiv A_{1}BD &\mathrm{mod}\ 3,\\
		A_{1}D\equiv 0 &\mathrm{mod}\ 4.
	\end{dcases}
\end{equation}
Next we compute all possible values of $s_{1}\left(N_{e}\right)$ when $N$ is nonspin. In formula \eqref{eq:s^I inv of N_e, N nonspin} if we set $s_{1}\left(N_{e}\right)=\frac{r}{28}\ \mathrm{mod}\ 1$, then
\begin{align*} 			   	  r&=\left(-18\left(8C_{1}u^{2}-4Buv+2A_{1}v^{2}\right)+6\left(B^{2}+6BC_{1}+2BD\right)u-6\left(A_{1}B+6A_{1}C_{1}+2A_{1}D\right)v\right)\\
&\quad-A_{1}D\left(B+6C_{1}+D\right)-A_{1}C_{1}\left(A_{1}+3B+9C_{1}\right)\ \mathrm{mod}\ 28.
\end{align*}
In the expression of $r$, the first two brackets are even and the last one is odd if and only if $A_{1}$ and $C_{1}$ are both odd. Hence $r$ is odd if and only if $A_{1}$ and $C_{1}$ are both odd. Since $B$ is odd we must have $B^{2}\equiv1\ \mathrm{mod}\ 8$. Hence from $B^{2}=4A_{1}C_{1}+1$ we deduce that $A_{1}C_{1}$ is even, and $r$ must be even. 
Now we consider all possible even values of $r$. From $B^{2}=4A_{1}C_{1}+1$ we simply take $A_{1}=B_{1}=1$, $C_{1}=D=0$, and $r\equiv72uv-36v^{2}+6u-6v\ \mathrm{mod}\ 28.$ Thus $r\ \mathrm{mod}\ 4=2(u-v)\ \mathrm{mod}\ 4$ and $r\ \mathrm{mod}\ 7=u^{2}-(u-v-3)^{2}+2\ \mathrm{mod}\ 7$. It is straightforward to see that $r\ \mathrm{mod}\ 4$ takes all even values of $\mathbb{Z}/4$. By Lemma \ref{Lemma from number theory} it is also true that $r\ \mathrm{mod}\ 7$ takes all values of $\mathbb{Z}/7$ since $\pm1$ are coprime to $7$. Hence by Chinese remainder theorem $r\ \mathrm{mod}\ 28$ takes all even values in $\mathbb{Z}/28$, i.e. $\mathcal{R}^{-}=2\mathbb{Z}/28$. The first statement of Proposition \ref{Propotision: homeom S2 times S5} is proved.

When $N$ is spin and $\det M_{e}=-1$, we have $A$ is even and $B$ and $C$ are both odd with $AC-B^{2}=-1$. First we apply formula \eqref{eq:s^I inv of N_e, N spin and det M_e=-1 } to equation \eqref{Eq: which N_e is homeom S2 times S5}. The $2-$primary part gives $A\equiv0\ \mathrm{mod}\ 8$ and $D\equiv1\ \mathrm{mod}\ 2$. Hence we set $A=8A_{1}$ and $B^{2}=8A_{1}C+1$. The $3-$primary part gives $A_{1}\left(C^{2}-BD-D^{2}\right)\equiv0\ \mathrm{mod}\ 3$. Next we compute all possible values of $s_{1}\left(N_{e}\right)$. In formula \eqref{eq:s^I inv of N_e, N spin and det M_e=-1 } if we set $s_{1}\left(N_{e}\right)=\frac{r}{28}\ \mathrm{mod}\ 1$, then
\begin{align*}
	r & = -18\left(Cu^{2}-2Buv+8A_{1}v^{2}\right)-2(3B(B-D)-2)u+48A_{1}(B-D)v+A_{1}-4A_{1}(B-D)^{2}\ \mathrm{mod}\ 28.
\end{align*}
Hence $r\ \mathrm{mod}\ 4=2u^{2}+A_{1}\ \mathrm{mod}\ 4$. When $A_{1}$ is odd $r\ \mathrm{mod}\ 4$ tkes all odds modulo $4$ and the same statement holds for even case. We take $A_{1}=1$, $B=3$, $C=D=1$ when $A_{1}$ is odd, and $r\ \mathrm{mod}\ 7=-(2u+v+1)^{2}+2v^{2}\ \mathrm{mod}\ 7$. We take $A_{1}=0$, $B=C=D=1$ when $A_{1}$ is even, and $r\ \mathrm{mod}\ 7=-(2u+v-1)^{2}+(v-1)^{2}\ \mathrm{mod}\ 7$. By Lemma \ref{Lemma from number theory}, in either case $r$ can take all values in $\mathbb{Z}/7$. Hence $r\ \mathrm{mod}\ 28$ takes all even values in $\mathbb{Z}/28$, i.e. $\mathcal{R}^{+}\supset\mathbb{Z}/28$. While $\mathcal{R}^{+}\subset\mathbb{Z}/28$, we obtain $\mathcal{R}^{+}=\mathbb{Z}/28$. This proves the second statement of Proposition \ref{Propotision: homeom S2 times S5} and hence we complete the proof of Proposition \ref{Propotision: homeom S2 times S5}.
\end{pf}
\begin{pf}[of Lemma \ref{Lemma from number theory}]
	Note that $(x,y)$ is a solution to the equation $ax^{2}+by^{2}\equiv c\ \mathrm{mod}\ p$ if and only if $(x,y)$ is a solution to the equation set $x^{2}\equiv u$, $y^{2}\equiv b^{-1}\left(c-au\right)\ \mathrm{mod}\ p$ for some $u$. It follows from the definition of Legendre symbol that the number of solutions to the equation $x^{2}\equiv u\ \mathrm{mod}\ p$ is $\left(\frac{u}{p}\right)+1$. Hence
	\begin{equation*}
		\begin{aligned}
			N(a,b,c;p)&=\sum_{u=0}^{p-1}\left(\left(\frac{u}{p}\right)+1\right)\left(\left(\frac{b^{-1}(c-au)}{p}\right)+1\right)\\
			&=\sum_{u=0}^{p-1}\left(\frac{u}{p}\right)\left(\frac{b^{-1}(c-au)}{p}\right)+\sum_{u=0}^{p-1}\left(\frac{u}{p}\right)+\sum_{u=0}^{p-1}\left(\frac{b^{-1}(c-au)}{p}\right)+p.
		\end{aligned}
	\end{equation*}
	We have $\sum_{u=0}^{p-1}\left(\frac{u}{p}\right)=\sum_{u=0}^{p-1}\left(\frac{b^{-1}(c-au)}{p}\right)=0$. Meanwhile when $u$ is coprime to $p$ we have $$\left(\frac{u}{p}\right)\left(\frac{b^{-1}(c-au)}{p}\right)=\left(\frac{b}{p}\right)\left(\frac{cu^{-1}-a}{p}\right).$$
	Hence if $c$ is divisible by $p$ we have $N(a,b,c;p)=(p-1)\left(\frac{b}{p}\right)\left(\frac{-a}{p}\right)+p=p+(p-1)\left(\frac{-ab}{p}\right)$, and if $c$ is coprime to $p$ we have $N(a,b,c;p)=\left(\frac{b}{p}\right)\left(-\left(\frac{-a}{p}\right)\right)+p=p-\left(\frac{-ab}{p}\right)$. Since the Legendre symbol takes value in $\pm1$ and $0$, the number $N(a,b,c;p)$ is always positive.
\end{pf}

\section{Circle actions on $S^{2}\times S^{5}\# S^{3}\times S^{4}\#\Sigma$}\label{Homeom S2 times S5 sum S3 times S4}

In this section we determine for which homotopy $7-$sphere $\Sigma$ the manifold $S^{2}\times S^{5}\#S^{3}\times S^{4}\#\Sigma$ admits a free circle action. The organization of this section is similar to that of Section \ref{Homeom S2 times S5}.

First in Section \ref{Subsec: Orbit B 2} we deal with the broader class of manifolds, i.e. the string cohomology $S^{2}\times S^{5}\#S^{3}\times S^{4}$. In particular when $\Sigma$ is a homotopy $7-$sphere the manifold $S^{2}\times S^{5}\#S^{3}\times S^{4}\#\Sigma$ is a string cohomology $S^{2}\times S^{5}\#S^{3}\times S^{4}$. We will study the invariants of the orbit space $N$ of some free circle action on some string cohomology $S^{2}\times S^{5}\#S^{3}\times S^{4}$. Next in Section \ref{Subsec: s inv 2} we compute the diffeomorphism invariants of those string cohomology $S^{2}\times S^{5}\#S^{3}\times S^{4}$ that admits a free circle action. We express the diffeomorphism invariants of such a manifold $M$ in terms of the invariants of the orbit space $N$. Finally in Section \ref{Subsec: compare s inv 2} we compare the diffeomorphism invariants of $M$ and $S^{2}\times S^{5}\#S^{3}\times S^{4}\#\Sigma$, and we determine for which homotopy $7-$sphere $\Sigma$ the manifold $S^{2}\times S^{5}\#S^{3}\times S^{4}\#\Sigma$ admits a free circle action.

Here a string manifold is defined as follows. When $M$ is a spin manifold, its \emph{first spin Pontryagin class} $\overline{p_{1}}(M)\in H^{4}(M)$ is defined and is related to the first Pontryagin class by $2\overline{p_{1}}(M)=p_{1}(M)$. We say $M$ is \emph{string} if it is spin and $\overline{p_{1}}(M)=0$. It is useful to notice that when $M$ is spin and its cohomology groups are torsion-free, the class $\overline{p_{1}}(M)$ vanishes if and only if $p_{1}(M)$ does and $M$ is string if and only if $p_{1}(M)=0$. 

%Now assume $M$ is a cohomology $S^{2}\times S^{5}\# S^{3}\times S^{4}$ that admits a free circle action. The organization of this section is the same as Section \ref{Homeom S2 times S5}. First we study the invariants of the orbit $N$. Then we express the $s-$invariant of $M$ (\cite[Section 4]{Kreck2018OnTC}) in terms of invariants of $N$. Finally we compare the invariants of $M$ and $S^{2}\times S^{5}\# S^{3}\times S^{4}\#\Sigma$ and determine for which $\Sigma\in\Theta_{7}$ the manifold $S^{2}\times S^{5}\# S^{3}\times S^{4}\#\Sigma$ admits a free circle action. For convenience a string cohomology $S^{2}\times S^{5}\# S^{3}\times S^{4}$ is called a \emph{manifold of type II}.

\subsection{The orbit space $N$}\label{Subsec: Orbit B 2}
In this section we study the invariants of the orbit space of a free circle action on $M$, where $M$ is a string cohomology $S^{2}\times S^{5}\#S^{3}\times S^{4}$. Again by \cite[Chapter 21, Theorem 21.10]{JohnLee} the orbit is a $6-$manifold and the projection $M\xrightarrow{p}N$ is a circle bundle. 

\begin{lem}\label{Lemma: Base space 2}
	Assume the $6-$manifold $N$ is the orbit of a free circle action on a string cohomology $S^{2}\times S^{5}\# S^{3}\times S^{4}$. Then $N$ is simply connected and satisfies the following conditions.
	\begin{compactenum}
		\item Its cohomology groups are torsion-free and the only nonzero Betti numbers are $b_{0}=b_{6}=1,\ b_{2}=b_{4}=2$;
		\item Its second cohomology group has a basis $(e,f)$ such that $\mu(e,e,e)=\varepsilon$, $\mu(e,e,f)=A$, $\mu(e,f,f)=\varepsilon A^2$ and $\mu(f,f,f)=A^{3}$, where $\varepsilon=\pm1$.
		\item It is spin and $p_{1}(N)=(4\varepsilon+24u)\left(e\spcheck+\varepsilon Af\spcheck\right)$ for some $u\in\mathbb{Z}$, where $\left(e\spcheck,f\spcheck\right)$ is the dual basis of $H^{4}(N)$ with respect to $(e,f)$.
	\end{compactenum}
\end{lem}
\begin{pf}
Let $M$ be a string cohomology $S^{2}\times S^{5}\#S^{3}\times S^{4}$ that admits a free circle action and let $N$ be the orbit manifold. Then we have a circle bundle $M\xrightarrow{p}N$. First we determine the cohomology groups of $N$. The long exact sequence of homotopy groups associated to this fiber bundle shows that $N$ is simply connected and there is a short exact sequence $0\to\pi_{2}(M)\to\pi_{2}(N)\to\pi_{1}\left(S^{1}\right)\to0$. Combining $\pi_{2}(M)\cong H_{2}(M)\cong\mathbb{Z}$ and $\pi_{1}\left(S^{1}\right)\cong\mathbb{Z}$ we have $\pi_{2}(N)\cong\mathbb{Z}^{2}$ and $H_{2}(N)\cong H^{2}(N)\cong\mathbb{Z}^{2}$. Let $e\in H^{2}(N)$ denote the Euler class of the circle bundle. From Gysin sequence we deduce the following two exact sequences:
\begin{align}
	0 \to H^{0}(N) &\xrightarrow{\cup e} H^{2}(N) \xrightarrow{p^{*}} H^{2}(M) \to 0,\label{SES: e primitive, again}\\
	0\to H^{3}(N)\xrightarrow{p^{*}}H^{3}(M)\to H^{2}(N)&\xrightarrow{\cup e} H^{4}(N) \xrightarrow{p^{*}} H^{4}(M)\to H^{3}(N)\to0.\label{LES: cup e on H2}
\end{align}
It follows from exact sequence \eqref{SES: e primitive, again} that $e$ is primitive. From exact sequence \eqref{LES: cup e on H2} we see that $H^{3}(N)$ injects into $H^{3}(M)\cong\mathbb{Z}$ and thus must be $0$ or $\mathbb{Z}$. Now $N$ has torsion-free cohomology groups. The rank of $H^{3}(N)$ is even by \cite[Theorem 1]{Wall66}, hence $H^{3}(N)=0$. Now we obtain that the only nonzero Betti numbers are $b_{0}=b_{6}=1$, $b_{2}=b_{4}=2$. This justifies Condition 1 of Lemma \ref{Lemma: Base space 2}.

Next we determine the cohomology ring of $N$. Back to exact sequence \eqref{LES: cup e on H2}, and we see the morphism $H^{2}(N)\xrightarrow{\cup e}H^{4}(N)$ has $\ker\cong H^{3}(M)\cong\mathbb{Z}$ and $\mathrm{coker}\cong H^{4}(M)\cong\mathbb{Z}$. Hence its image in $H^{4}(N)$ is a direct summand and is isomorphic to $\mathbb{Z}$. We extend $e\in H^{2}(N)$ to a basis $(e,f)$ and equip $H^{4}(N)$ with the dual basis $\left(e\spcheck,f\spcheck\right)$. With respect to these bases the morphism $H^{2}(N)\xrightarrow{\cup e}H^{4}(N)$ is represented by the matrix $\ M_{e}=\begin{pmatrix}
	\mu(e,e,e) & \mu(e,f,e) \\ 
	\mu(e,e,f) & \mu(e,f,f)
\end{pmatrix}$. Since the image of $H^{2}(N)\xrightarrow{\cup e}H^{4}(N)$ is a direct summand isomorphic to $\mathbb{Z}$, the products $e^{2}$ and $ef$ are colinear. Hence $\det M_{e}=\mu(e,e,f)^{2}-\mu(e,e,e)\mu(e,f,f)=0$ and at least one column of $M_{e}$ is primitive in $\mathbb{Z}^{2}$. This implies at least one of the greatest common divisors $\mathrm{gcd}\left(\mu(e,e,e),\mu(e,e,f)\right)$ and $\mathrm{gcd}\left(\mu(e,e,f),\mu(e,f,f)\right)$ is $1$. Hence we obtain that either $\mu(e,e,e)=\pm1$, $\mu(e,f,f)=\mu(e,e,e)\mu(e,e,f)^{2}$ or $\mu(e,f,f)=\pm1$, $\mu(e,e,e)=\mu(e,f,f)\mu(e,e,f)^{2}$. The product $f^{2}=\mu(f,f,e)e\spcheck+\mu(f,f,f)f\spcheck$ is not contained above. Since $e$ is primitive, $p^{*}f$ generates $H^{2}(M)\cong\mathbb{Z}$. Since $M$ is a cohomology $S^{2}\times S^{5}\#S^{3}\times S^{4}$ we have $p^{*}\left(f^{2}\right)=\left(p^{*}f\right)^{2}=0$. Hence by exact sequence \eqref{LES: cup e on H2} we see $f^{2}\in\mathrm{im}\left(H^{2}(N)\xrightarrow{\cup e}H^{4}(N)\right)$, and the products $e^{2},\ ef,\ f^{2}$ are all colinear. Write
\begin{equation*}
	\begin{pmatrix}
		e^{2} & ef & f^{2}
	\end{pmatrix}=\begin{pmatrix}
	e\spcheck & f\spcheck
	\end{pmatrix}\begin{pmatrix}
	\mu(e,e,e) & \mu(e,e,f) & \mu(e,f,f) \\
	\mu(e,e,f) & \mu(e,f,f) & \mu(f,f,f)
	\end{pmatrix},
\end{equation*}
and the $2\times3$ matrix has rank $1$. It is routine to deduce if $\mu(e,f,f)=\pm1$, $\mu(e,e,e)$, $\mu(e,e,f)$ and $\mu(f,f,f)$ are all $\pm1$. Hence we set $\varepsilon=\mu(e,e,e)=\pm1$ and $A=\mu(e,e,f)\in\mathbb{Z}$, then $ef=Ae^{2}$, $f^{2}=A^{2}e^{2}$ and $\mu(e,f,f)=\varepsilon A^{2}$, $\mu(f,f,f)=A^{3}$. This justifies Condition 2 of Lemma \ref{Lemma: Base space 2}.

Now we study the characteristic classes of $N$. We keep using diagram \eqref{BundleDiagram} and the previous notations. As before we begin with $w_{2}(N)$. From $w_{2}\left(M\right)=0$ we still have $w_{2}(N)=0$ or $\overline{e}$. If $w_{2}(N)=\overline{e}$, formula \eqref{Formula: Wu formula} would imply $0=\overline{\mu(e,e,e)}=1\ \mathrm{mod}\ 2$. This is a contradiction. Hence $N$ is spin. Next we determine $p_{1}(N)$. Set $p_{1}(N)=k\cdot e\spcheck+l\cdot f\spcheck$, $k,l\in\mathbb{Z}$. Since $M$ is string, we have $p_{1}\left(M\right)=p^{*}\left(p_{1}(N)+e^{2}\right)=p^{*}\left(p_{1}(N)\right)=0$. By exact sequence \eqref{LES: cup e on H2} we see $p_{1}(N)\in\mathrm{im}\left(H^{2}(N)\xrightarrow{\cup e}H^{4}(N)\right)$. Since $e^{2}=\varepsilon e\spcheck+A f\spcheck$ and $ef=A e\spcheck+\varepsilon A^{2}f\spcheck=\varepsilon A e^{2}$, the image of $H^{2}(N)\xrightarrow{\cup e}H^{4}(N)$ is generated by $e^{2}$. Hence $(k,l)$ and $(\varepsilon,A)$ are colinear. Now we conclude that $l=\varepsilon Ak$ and $p_{1}(N)=k\varepsilon e^{2}=k\left(e\spcheck+\varepsilon Af\spcheck\right)$. Combining formula \eqref{Formula: p1 and w2} and $w_{2}(N)=0$ we see $k$ is even. Now we unpack formula \eqref{Jupp relation}. As $N$ is spin, we take $\widehat{w}=0$. Set $x=Xe+Yf$, $X,Y\in\mathbb{Z}$, and we have
\begin{equation*}
	4X^{3}\varepsilon+12X^{2}YA+12XY^{2}\varepsilon A^{2}+4Y^{3}A^{3}\equiv kX+ \varepsilon AkY\ \mathrm{mod}\ 24,\ \forall X,Y\in\mathbb{Z}.
\end{equation*}
Solving this equation we obtain $k\equiv4\varepsilon\ \mathrm{mod}\ 24$. This justifies Condition 3 of Lemma \ref{Lemma: Base space 2}. Now we complete the proof.
\end{pf}
\begin{lem}\label{Lemma: Base space 2, if part}
	Let $\varepsilon=\pm1$ and let $A$ be an integer. Up to orientation-preserving diffeomorphism there is a unique simply connected spin $6-$manifold $N$ that satisfies Conditions 1, 2 and 3 of Lemma \ref{Lemma: Base space 2} and is the orbit of a free circle action on a string cohomology $S^{2}\times S^{5}\#S^{3}\times S^{4}$.
\end{lem}
\begin{pf}
	The proof is parallel to the proof of Lemma \ref{Lemma: Base space 1, if part}. By \cite[Theorem 1]{Jupp1973} up to orientation-preserving diffeomorphism there is a unique simply connected spin $6-$manifold $N$ which satisfies Conditions 1, 2 and 3 of Lemma \ref{Lemma: Base space 2}. Let $N_{e}$ be the total space of the circle bundle over $N$ whose Euler class is $e$. Imitating the proof of Lemma \ref{Lemma: Base space 1, if part}, we see that $N_{e}$ is a string cohomology $S^{2}\times S^{5}\#S^{3}\times S^{4}$.
\end{pf}

\subsection{Invariants of $N_{e}$}\label{Subsec: s inv 2}
In this section we compute the diffeomorphism invariants of $N_{e}$. The diffeomorphism invariants are the $s-$invariants developed by Kreck (\cite[Section 4]{Kreck2018OnTC}). First we review the definition of the original $s-$invariants. Then we extend the definition. Finally we apply the extended definition of $s-$invariants to the manifold $N_{e}$.

\cite[Section 4]{Kreck2018OnTC} considered the diffeomorphism invariant of spin cohomology $S^{2}\times S^{5}\#S^{3}\times S^{4}$. Here for convenience a string cohomology $S^{2}\times S^{5}\#S^{3}\times S^{4}$ is called a \emph{manifold of type II}. Then given a pair $(N,e)$ satisfying Lemma \ref{Lemma: Base space 2}, the manifold $N_{e}$ is of type II. A generator of $H^{2}\left(N_{e}\right)\cong\mathbb{Z}$ is $p^{*}f$. When $\Sigma$ is a homotopy $7-$sphere, the manifold $S^{2}\times S^{5}\#S^{3}\times S^{4}\#\Sigma$ is also a type II manifold. 

The diffeomorphism invariants of a type II manifold are defined as follows. Let $M$ be a type II manifold and let $z$ be a generator of $H^{2}(M)$. By \cite[Theorem 6]{Kreck2018OnTC} $M$ admits a spin coboundary $W$ such that there is a class $\widehat{z}\in H^{2}(W)$ with $\widehat{z}|_{M}=z$. We have the following characteristic numbers
\begin{eqnarray}\label{Formula: inv of type II manifolds, spin coboundary}
	\left\{\begin{aligned}
		S_{1}\left(W,\widehat{z}\right)&=\left<\overline{p_{1}}^{2},[W,M]\right>\in\mathbb{Z},\\
		S_{2}\left(W,\widehat{z}\right)&=\left<\widehat{z}^{2}\left(\widehat{z}^{2}-\overline{p_{1}}\right),[W,M]\right>\in\mathbb{Z},\\
		S_{3}\left(W,\widehat{z}\right)&=\left<\widehat{z}^{2}\overline{p_{1}},[W,M]\right>\in\mathbb{Z}.
	\end{aligned}\right.
\end{eqnarray}
Here $\overline{p_{1}}=\overline{p_{1}}(W)$. The integral classes $\overline{p_{1}}$ and $\widehat{z}^{2}$ vanish when they restrict to $M$, hence they can be lifted to $H^{4}(W,M)$ and the products $\overline{p_{1}}^{2}$, $\widehat{z}^{2}\overline{p_{1}}$ and $\widehat{z}^{4}$ are classes in $H^{8}(W,M)$ evaluating on $[W,M]$. For convenience we abbreviate $\left<\overline{p_{1}}^{2},[W,M]\right>$ for $\overline{p_{1}}^{2}$ and the other monomial terms are silimar. Then we have
\begin{eqnarray}\label{Formula: inv of type II manifolds, spin coboundary, expand}
	\left\{\begin{aligned}
		S_{1}\left(W,\widehat{z}\right)&=\overline{p_{1}}^{2}\in\mathbb{Z},\\
		S_{2}\left(W,\widehat{z}\right)&=\widehat{z}^{4}-\widehat{z}^{2}\overline{p_{1}}\in\mathbb{Z},\\
		S_{3}\left(W,\widehat{z}\right)&=\widehat{z}^{2}\overline{p_{1}}\in\mathbb{Z}.
	\end{aligned}\right.
\end{eqnarray}
Set 
\begin{alignat*}{7}
		s_{1}(M)&=&S_{1}\left(W,\widehat{z}\right)&\ \mathrm{mod}&\ 224&\in&\ \mathbb{Z}/224,\\
		s_{2}(M)&=&S_{2}\left(W,\widehat{z}\right)&\ \mathrm{mod}&\ 24&\in&\ \mathbb{Z}/24,\\
		s_{3}(M)&=&S_{2}\left(W,\widehat{z}\right)&\ \mathrm{mod}&\ 2&\in&\ \mathbb{Z}/2,
	\end{alignat*}
and these invariants $s_{i}$ do not depend on the choices of the generators of $H^{2}(M)$ or the spin coboundaries $\left(W,\widehat{z}\right)$. They are invariants of the type II manifold $M$. It follows from \cite[Theorem 2]{Kreck2018OnTC} that two manifolds $M$, $M'$ of type II are diffeomorphic if and only if they share the same invariants $s_{1}$, $s_{2}$ and $s_{3}$.

Sometimes we cannot easily construct such a spin coboundary, but we can find a nonspin coboundary $W$ and its cohomology classes $\widehat{z},c\in H^{2}(W)$ such that $w_{2}(W)=\overline{c}$ and $\left.\widehat{z}\right|_{M}=z,\ c|_{M}=0$.
Moreover we can still use the tripple $\left(W,\widehat{z},c\right)$ to compute the $s-$invariants of $M$. Hence we first extend the previous definition (\cite[Section 4]{Kreck2018OnTC}) of $s-$invariants for type II manifolds. We rewrite the formulae of $s-$invariants of type II manifolds and add the superscript II.
\begin{alignat*}{7}
	s_{1}^{\mathrm{II}}(M) &=& \frac{1}{896}p_{1}^{2} &\in& \ \left.\frac{1}{896}\mathbb{Z}\right/\mathbb{Z} &\subset& \mathbb{Q}/\mathbb{Z}, \\
	s_{2}^{\mathrm{II}}(M) &=& \frac{1}{24}\widehat{z}^{4}-\frac{1}{48}\widehat{z}^{2}p_{1} &\in&  \ \left.\frac{1}{24}\mathbb{Z}\right/\mathbb{Z} &\subset& \mathbb{Q}/\mathbb{Z},\\
	s_{3}^{\mathrm{II}}(M) &=& \frac{1}{4}\widehat{z}^{2}p_{1} &\in &\left.\frac{1}{2}\mathbb{Z}\right/\mathbb{Z} &\subset & \mathbb{Q}/\mathbb{Z}.
\end{alignat*}
Here we use the first Pontryagin class since the first spin Pontryagin class is not defined for nonspin coboundary. Also we recall the expression of $s-$invariants for Type I manifolds. Here we add the superscript I and assume the coboundary has signature zero.
\begin{alignat*}{6}
	s^{\mathrm{I}}_{1}(M) &\  = &\  &\frac{1}{896}p_{1}^{2} &\  \in &\  \mathbb{Q}/\mathbb{Z},\\
	s^{\mathrm{I}}_{2}(M) &\  = &\  &\frac{1}{24}\widehat{z}^{4}-\frac{1}{48}\widehat{z}^{2}p_{1} &\  \in &\  \mathbb{Q}/\mathbb{Z},\\
	s^{\mathrm{I}}_{3}(M) &\  = &\  &\frac{2}{3}\widehat{z}^{4}-\frac{1}{12}\widehat{z}^{2}p_{1} &\  \in &\  \mathbb{Q}/\mathbb{Z}.
\end{alignat*}
The invariants $s^{\mathrm{I}}$ and $s^{\mathrm{II}}$ are rational linear combinations of $p_{1}^{2},\ \widehat{z}^{2}p_{1}$ and $\widehat{z}^{4}$ such that
\begin{equation*}
	\begin{pmatrix}
		s_{1}^{\mathrm{II}} & s_{2}^{\mathrm{II}} & s_{3}^{\mathrm{II}}
	\end{pmatrix}=\begin{pmatrix}
		s_{1}^{\mathrm{I}} & s_{2}^{\mathrm{I}} & s_{3}^{\mathrm{I}}
	\end{pmatrix}\begin{pmatrix}
		1 & 0 & 0\\
		0 & 1 & 0\\
		0 & -16 & 1
	\end{pmatrix}.
\end{equation*}
The transition matrix is lower triangular and its diagonal entries are all $1$. Hence this matrix is invertible over $\mathbb{Z}$, and two type II manifolds have the same $s^{\mathrm{II}}-$invariants if and only if they have the same $s^{\mathrm{I}}-$invariants. We conclude that if the coboundaries are spin and have signature zero, the diffeomorphism invariants for type I manifolds and type II manifolds have the same expression. Hence we drop the superscripts $\mathrm{I},\ \mathrm{II}$ and write $s_{i}=s_{i}^{\mathrm{I}}=s_{i}^{\mathrm{II}}$ uniformly. Further by \cite[Remark 2.6]{Kreck1991SomeNH}, given a type II manifold $M$ with $z$ a generator of $H^{2}(M)$, if $W$ is a nonspin coboundary of type II manifold $M$ and there are classes $\widehat{z},\ c\in H^{2}(W)$ such that $w_{2}(W)=\overline{c}$, $\widehat{z}|_{M}=z$ and $c|_{M}=0$, then formula \eqref{Formula: s invariant, expand} still apply to the computation of $s-$invariants of $M$.

\begin{eg}
	Let $(N,e)$ be given as in Lemma \ref{Lemma: Base space 2}. The $s-$invariants of the type II manifold $N_{e}$ are given by:
	\begin{eqnarray}\label{eq:s^II inv of N_e}
		\left\{\begin{aligned}
			s_{1}\left(N_{e}\right) & = & \frac{9\varepsilon u^{2}+2u}{14}\ \mathrm{mod}\ 1,\\
			s_{2}\left(N_{e}\right) & = & 0\ \mathrm{mod}\ 1,\\
			s_{3}\left(N_{e}\right) & = & 0\ \mathrm{mod}\ 1.
		\end{aligned}
		\right.
	\end{eqnarray}
	
The $s-$invariants of $N_{e}$ is computed as follows. A natural coboundary of $N_{e}$ is the associated disc bundle $D_{e}$. Note that $D_{e}$ is nonspin and $w_{2}\left(N_{e}\right)=\overline{\pi^{*}e}$. Set $z=p^{*}f$, $\widehat{z}=\pi^{*}f$, $c=\pi^{*}e$, and $\left(D_{e},\widehat{z},c\right)$ is a nonspin coboundary of $(N_{e},z)$. By Lemma \ref{Lemma: Base space 2} we have $p_{1}(N)=(4+24\varepsilon u)e^2$, where $\varepsilon=\pm1$ and $u\in\mathbb{Z}$. Then $p_{1}\left(D_{e}\right)=\pi^{*}\left(p_{1}(N)+e^{2}\right)=(1+\varepsilon k)\pi^{*}\left(e^{2}\right)=(5+24\varepsilon u)\pi^{*}\left(e^{2}\right)$. Just as before the signature term is $\sigma\left(D_{e},N_{e}\right)=\sigma\left(M_{e}\right)=\sigma\begin{pmatrix}
\varepsilon & A \\ A & \varepsilon A^{2}
\end{pmatrix}=\varepsilon$. The remaining monomials are
\begin{alignat*}{5}
	p_{1}\cdot p_{1} & = & \mu\left(e,(5+24\varepsilon u)e,(5+24\varepsilon u)e\right) & = &(5+24\varepsilon u)^{2}\varepsilon;\\
	c^{2}\cdot p_{1} & = & \mu(e,e,(5+24\varepsilon u)e) & = & (5+24\varepsilon u)\varepsilon,\\
	\widehat{z}c\cdot p_{1} & = & \mu(e,f,(5+24\varepsilon u)e) & = & (5+24\varepsilon u)A,\\
	\widehat{z}^{2}\cdot p_{1} & = & \mu\left(e,A^{2}e,(5+24\varepsilon u)e\right) & = & (5+24\varepsilon u)\varepsilon A^{2};\\
	c^{2}\cdot c^{2} & = & \mu(e,e,e) & = & \varepsilon,\\
	\widehat{z}c\cdot c^{2} & = & \mu(e,f,e) & = & A,\\
	\widehat{z}^{2}\cdot c^{2} & = & \mu(e,A^{2}e,e) & = & \varepsilon A^{2},\\
	\widehat{z}c\cdot \widehat{z}c & = & \mu(e,f,f) & = & \varepsilon A^{2},\\
	\widehat{z}^{2}\cdot \widehat{z}c & = & \mu\left(e,A^{2}e,f\right) & = & A^{3},\\
	\widehat{z}^{2}\cdot \widehat{z}^{2} & = & \mu\left(e,A^{2}e,A^{2}e\right) & = & \varepsilon A^{4}.
\end{alignat*}
Here we have two ways to compute the monomial $\widehat{z}^{2}c^{2}$, and they lead the same result. Hence
\begin{eqnarray*}
	S_{1}\left(D_{e},\widehat{z},c\right) & = & -\frac{1}{224}\varepsilon+\frac{1}{896}(5+24\varepsilon u)^{2}\varepsilon-\frac{1}{192}(5+24\varepsilon u)\varepsilon+\frac{1}{384}\varepsilon\\
	& = & \frac{1}{7}u+\frac{9}{14}\varepsilon u^{2},\\
	S_{2}\left(D_{e},\widehat{z},c\right) & = &-\frac{1}{48}\left((5+24\varepsilon u)\varepsilon A^{2}+(5+24\varepsilon u)A\right)+\frac{1}{48}\left(2\varepsilon A^{4}+4 A^{3}+3\varepsilon A^{2}+A\right)\\
	& = & \frac{1}{24}\varepsilon A(A+1)(A-1)(A+2\varepsilon)-\frac{1}{2}uA(A+\varepsilon)\in\mathbb{Z},\\
	S_{3}\left(D_{e},\widehat{z},c\right) & = &-\frac{1}{24}\left(2(5+24\varepsilon u)\varepsilon A^{2}+ (5+24\varepsilon u)A\right)+\frac{1}{24}\left(16\varepsilon A^{4}+16A^{3}+6\varepsilon A^{2}+A\right)\\
	& = & \frac{1}{6}\varepsilon A(A+\varepsilon)(2A+1)(2A-1)-\varepsilon uA(2\varepsilon A+1)\in\mathbb{Z}.
\end{eqnarray*}
Modulo $\mathbb{Z}$ and we obtain formula \eqref{eq:s^II inv of N_e}.

It follows from \cite[Section 3]{MontgomeryYang68} that the homotopy $7-$sphere $\Sigma$ admits a free circle action if and only if $\mu(\Sigma)\in\left\{ \frac{9\varepsilon u^{2}+2u}{14}\ \mathrm{mod}\ 1:\varepsilon=\pm1,u\in\mathbb{Z}\right\}.$ See also Remark \ref{Rmk: homotopy spheres admitting circle actions}.
\end{eg}

\subsection{Comparing $s-$invariants}\label{Subsec: compare s inv 2}
In this section we first compute the $s-$invariants of $S^{2}\times S^{5}\#S^{3}\times S^{4}\#\Sigma_{r}$, then we compare the $s-$invariants of the manifold $N_{e}$ of type II and the manifold $S^{2}\times S^{5}\#S^{3}\times S^{4}\#\Sigma_{r}$, deducing for which homotopy $7-$sphere the manifold $S^{2}\times S^{5}\#S^{3}\times S^{4}\#\Sigma$ admits a free circle action.

%From \cite[Section 3]{EKinv} we see the $\mu-$invariant is defined for $S^{3}\times S^{4}$ and we can use the coboundary $S^{3}\times D^{5}$ to compute that $\mu\left(S^{3}\times S^{4}\right)=0\ \mathrm{mod}\ 1$. Hence the $\mu-$invariant is also defined for the connected sum $S^{2}\times S^{5}\#S^{3}\times S^{4}\#\Sigma_{r}$ with $\mu\left(S^{2}\times S^{5}\#S^{3}\times S^{4}\#\Sigma_{r}\right)=\mu\left(S^{2}\times S^{5}\right)+\mu\left(S^{3}\times S^{4}\right)+\mu\left(\Sigma_{r}\right)=\frac{r}{28}\ \mathrm{mod}\ 1$. Note that for type II manifolds, the $s_{1}-$invariant also extends the $\mu-$invariant. This follows immediately from the expressions of invariants $\mu$ and $s_{1}$. Hence $s_{1}\left(S^{2}\times S^{5}\#S^{3}\times S^{4}\#\Sigma_{r}\right)=\mu\left(S^{2}\times S^{5}\#S^{3}\times S^{4}\#\Sigma_{r}\right)=\frac{r}{28}\ \mathrm{mod}\ 1$.

The boundary connected sum $S^{2}\times D^{6}\natural D^{4}\times S^{4}\natural W_{r}$ is a coboundary of $S^{2}\times S^{5}\# S^{3}\times S^{4}\#\Sigma_{r}$. By definition we can use the coboundary $S^{2}\times D^{6}\natural D^{4}\times S^{4}\natural W_{r}$ tocompute the $s-$invariants of $S^{2}\times S^{5}\# S^{3}\times S^{4}\#\Sigma_{r}$ and the result is as follows.
	\begin{eqnarray}\label{eq:s^II inv of homeom S2 times S5 sum S3 times S4}
	\left\{\begin{aligned}
	s_{1}\left(S^{2}\times S^{5}\# S^{3}\times S^{4}\#\Sigma_{r}\right)&=&\frac{r}{28}\ \mathrm{mod}\ 1,\\
	s_{2}\left(S^{2}\times S^{5}\# S^{3}\times S^{4}\#\Sigma_{r}\right)&=&0\ \mathrm{mod}\ 1,\\
	s_{3}\left(S^{2}\times S^{5}\# S^{3}\times S^{4}\#\Sigma_{r}\right)&=&0\ \mathrm{mod}\ 1.
	\end{aligned}
	\right.
\end{eqnarray}
Note that by the original definition for a type II manifold $M$ its $\mu$-invariant and $s_{1}-$invariant cannot be defined simultaneously. By \cite[Section 3, Condition $\mu$(a)]{EKinv} the $\mu-$invariant is defined for $M$ if there is a coboundary $W$ such that $H^{4}(W,M;\mathbb{Q})\to H^{4}(W;\mathbb{Q})$ is an isomorphism. While the $s_{1}-$invariant is defined for $M$ if there is a coboundary $W'$ such that $H^{4}(W')\to H^{4}(M)$ is an epimorphism.

Now we compare formulae \eqref{eq:s^II inv of N_e} and \eqref{eq:s^II inv of homeom S2 times S5 sum S3 times S4}. Immediately we obtain the following.
\begin{prop}\label{Proposition: homeom S2 times S5 sum S3 tims S4}
	Let $\Sigma$ be a homotopy $7-$sphere. The manifold $S^{2}\times S^{5}\#S^{3}\times S^{4}\#\Sigma$ admits a free circle action if and only if $\Sigma$ does. In this case the orbit must be spin. Furthermore, if a simply connected string cohomology $S^{2}\times S^{5}\#S^{3}\times S^{4}$ admits a free circle action, it must be diffeomorphic to $S^{2}\times S^{5}\#S^{3}\times S^{4}\#\Sigma$ for some homotopy $7-$sphere $\Sigma$ which admits a free circle action.
\end{prop}

\section{Circle actions on general $kS^{2}\times S^{5}\#lS^{3}\times S^{4}$}\label{More cnt sum}
%\cite[Theorem E]{galazgarcía2023free} implies that for any nonnegative integers $k,\ l$, the manifold $kS^{2}\times S^{5}\#lS^{3}\times S^{4}$ always admits a free circle action. 
In this section we consider the remaining cases. We will study when $k>1$ or $l>1$, for which homotopy $7-$sphere $\Sigma$ the manifold $kS^{2}\times S^{5}\#lS^{3}\times S^{4}\#\Sigma$ admits a free circle action. We will first recall some basic constructions about suspension operation from \cite{duan2022circle} and \cite{galazgarcía2023free}. Then we treat the case $k>1$ or $l>1$. Finally we combine all our results and complete the proof of our main theorem.

First we give some basic constructions. Given an $n-$manifold $N$, there are two associated $(n+1)-$manifolds $\Sigma_{0}N$ and $\Sigma_{1}N$ (\cite[Definition 1.1]{duan2022circle}). We denote $\Sigma N=\Sigma_{0}N=\Sigma_{1}N$ when $\Sigma_{0}N$ and $\Sigma_{1}N$ are diffeomorphic. By \cite[Proposition 3.2]{duan2022circle} we have $\Sigma\left(S^{2}\times S^{4}\right)=S^{2}\times S^{5}\#S^{3}\times S^{4}$ and $\Sigma\left(S^{3}\times S^{3}\right)=2S^{3}\times S^{4}$. If we take $e\in H^{2}(N)$, there is an associated $(n+1)-$manifold $\widetilde{\Sigma}_{e}N$ (\cite[Definition 5.1]{galazgarcía2023free}). By \cite[Theorem B(1)]{galazgarcía2023free} we have $\widetilde{\Sigma}_{h}\mathbb{C}P^{3}=S^{2}\times S^{5}$, where $h\in H^{2}\left(\mathrm{C}P^{3}\right)\cong\mathbb{Z}$ is the Euler class of Hopf bundle $S^{7}\to\mathbb{C}P^{3}$. For the connected sum of $6-$manifolds $N^{6}=N_{1}^{6}\#\cdots\#N_{m}^{6}$ we have $H^{2}\left(N^{6}\right)\cong\bigoplus_{i=1}^{m}H^{2}\left(N_{i}^{6}\right)$. When $\mathbb{C}P^{3}$ is the $i-$th summand, the collapsing map $N\to \mathbb{C}P^{3}$ pulls back $h\in H^{2}\left(\mathbb{C}P^{3}\right)$ to $h_{i}\in H^{2}(N)$.

Now we are prepared to study when $k>1$ or $l>1$, for which homotopy $7-$sphere $\Sigma$ the manifold $kS^{2}\times S^{5}\#lS^{3}\times S^{4}\#\Sigma$ admits a free circle action.

\begin{prop}\label{Prop: general}
	Let $\Sigma$ be a homotopy $7-$sphere.
	\begin{compactenum}
		\item If $k\geqslant2$, for any $l\in\mathbb{N}$ the manifold $kS^{2}\times S^{5}\#lS^{3}\times S^{4}\#\Sigma$ always admits a free circle action.
		\item If $l\geqslant2$ is even, the manifold $S^{2}\times S^{5}\#lS^{3}\times S^{4}\#\Sigma$ always admits a free circle action.
		\item If $l\geqslant2$ is odd, the manifold $S^{2}\times S^{5}\#lS^{3}\times S^{4}\#\Sigma$ admits a free circle action if and only if $\Sigma$ does.
	\end{compactenum}
\end{prop}
In the proof we will use the following lemma, whose proof is postponed later.
\begin{lem}\label{Lemma: k>1, direct construction}
	Let $k\geqslant2$, $l\geqslant0$, homotopy $7-$sphere $\Sigma$, $N_{\Sigma}$ and $e\in H^{2}\left(N_{\Sigma}\right)$ be given as above.
	\begin{compactenum}
		\item If $l$ is even, take $N=N_{\Sigma}\#(k-1)\mathbb{C}P^{3}\#\frac{l}{2}S^{3}\times S^{3}$, $e=e_{0}+h_{1}+\cdots+h_{k-1}\in H^{2}(N)$, and the total space of the circle bundle over $N$ with Euler class $e$ is $N_{e}=kS^{2}\times S^{5}\#lS^{3}\times S^{4}\#\Sigma$.
		\item If $l$ is odd, take $N=N_{\Sigma}\#(k-2)\mathbb{C}P^{3}\#S^{2}\times S^{4}\#\frac{l-1}{2}S^{3}\times S^{3}$, $e=e_{0}+h_{1}+\cdots+h_{k-2}\in H^{2}(N)$, and the total space of the circle bundle over $N$ with Euler class $e$ is $N_{e}=kS^{2}\times S^{5}\#lS^{3}\times S^{4}\#\Sigma$.
	\end{compactenum}
\end{lem}
\begin{pf}[of Proposition \ref{Prop: general}]
	First we assume $k\geqslant2$. By Proposition \ref{Propotision: homeom S2 times S5} the manifold $S^{2}\times S^{5}\#\Sigma$ admits a free circle action with a spin orbit. Let $N_{\Sigma}$ be such a spin orbit and let $e_{0}\in H^{2}\left(N_{\Sigma}\right)$ be the Euler class of associated circle bundle. In particular $\left(N_{\Sigma}\right)_{e_{0}}=S^{2}\times S^{5}\#\Sigma$. Then the proof of Proposition \ref{Prop: general}, statement 1 follows directly from Lemma \ref{Lemma: k>1, direct construction}.

	In the remaining we set $k=1$. When $l$ is even we give the construction directly. Let the homotopy $7-$sphere $\Sigma$, $N_{\Sigma}$, $e\in H^{2}\left(N_{\Sigma}\right)$ be as above. Set $N=N_{\Sigma}\#\frac{l}{2}S^{3}\times S^{3}$, and $H^{2}(N)\cong H^{2}\left(N_{\Sigma}\right)$. Identify $e\in H^{2}\left(N_{\Sigma}\right)$ with $e'\in H^{2}(N)$, and we have  $N_{e'}=\left(N_{\Sigma}\right)_{e}\#\Sigma\left(\frac{l}{2}S^{3}\times S^{3}\right)=S^{2}\times S^{5}\#\Sigma\#l S^{3}\times S^{4}$ by \cite[Theorem B, Proposition 3.2]{duan2022circle}. Hence for any even $l>0$ and any homotopy $7-$sphere $\Sigma$ the manfiold $S^{2}\times S^{5}\#lS^{3}\times S^{4}\#\Sigma$ admits a free circle action. This proves Proposition \ref{Prop: general}, statement 2.
	
	Now assume $l$ is odd. We can show as before that if the homotopy $7-$sphere $\Sigma$ admits a free circle action, so does $S^{2}\times S^{5}\#l S^{3}\times S^{4}\#\Sigma$. Take such homotopy $7-$sphere $\Sigma$. Let $N_{\Sigma}$ be an orbit of some free circle action on $\Sigma$, $e\in H^{2}\left(N_{\Sigma}\right)$ the Euler class and $\left(N_{\Sigma}\right)_{e}=\Sigma$. Set $N=N_{\Sigma}\#S^{2}\times S^{4}\#\frac{l-1}{2}S^{3}\times S^{3}$ and identify $e$ with $e'\in H^{2}(B)$. Then $N_{e'}=\left(N_{\Sigma}\right)_{e}\#\Sigma\left(S^{2}\times S^{4}\#\frac{l-1}{2}S^{3}\times S^{3}\right)=\Sigma\#S^{2}\times S^{5}\#lS^{3}\times S^{4}$ by \cite[Theorem B, Proposition 3.2]{duan2022circle}. The converse is also true. If $l$ is odd and $S^{2}\times S^{5}\#l S^{3}\times S^{4}\#\Sigma$ admits a free circle action, so does $\Sigma$. Suppose for a homotopy $7-$sphere $\Sigma$, the manifold $M=S^{2}\times S^{5}\#l S^{3}\times S^{4}\#\Sigma$ admits a free circle action. Let $N$ be such an orbit. Let $e\in H^{2}\left(N\right)$ be the Euler class of the associated circle bundle. Then $M=N_{e}$. The long exact sequence of homotopy groups associated to the circle bundle $M\xrightarrow{p}N$ implies that $N$ is simply connected. Hence $H^{1}\left(N\right)=H^{5}\left(N\right)=0$. The Gysin sequence associated to the circle bundle $M\xrightarrow{p}N$ implies that $N$ has torsion free cohomology groups and $H^{3}\left(N\right)$ embeds into $H^{3}\left(M\right)\cong\mathbb{Z}^{l}$. From \cite[Theorem 1]{Wall66} we see that $b_{3}\left(N\right)=2r$ is even and $N=N_{0}\#r S^{3}\times S^{3}$, where $N_{0}$ is another simply connected $6-$manifold with torsion free homology and vanishing third cohomology group. Let $N\xrightarrow{\rho}N_{0}$ be the collapsing map. Then $\rho$ induces isomorphism on the second cohomology groups, and there is a class $e_{0}\in H^{2}\left(N_{0}\right)$ such that $e=\rho^{*}e_{0}$. Also we set $M_{0}=\left(N_{0}\right)_{e_{0}}$ and let $M_{0}\xrightarrow{p_{0}}N_{0}$ denote the bundle projection. Then $M=M_{0}\#\Sigma\left(r S^{3}\times S^{3}\right)=M_{0}\#2r S^{3}\times S^{4}$, $b_{3}\left(M\right)=b_{3}\left(M_{0}\right)+2r=l$ and $b_{3}\left(M_{0}\right)$ is odd. From the Gysin sequence associated to the circle bundle $M_{0}\xrightarrow{p_{0}}N_{0}$. we obtain that $b_{3}\left(M_{0}\right)\leq b_{2}\left(N_{0}\right)=2$, hence $b_{3}\left(M_{0}\right)=1$ and $l=2r+1$. See the following diagram.
	$$\xymatrix{
		& H^{3}\left(N\right)\ar[r] & 0 & \\
		& H^{2}\left(N\right)\ar[r]^{\cup e} & H^{4}\left(N\right)\ar[r]^{p^{*}} & H^{4}\left(M\right)\ar@{->}`r[u] `u[l] `[lllu] `[ull] [ull] \\
		& 0 \ar[r]& H^{3}\left(N\right)\ar@{->}[r]^{p^{*}} & H^{3}\left(M\right)\ar`r[u] `u[l] `[lllu] `[ull] [ull] \\
		& H^{0}\left(N\right)\ar@{->}[r]^{\cup e} & H^{2}\left(N\right)\ar@{->}[r]^{p^{*}} & H^{2}\left(M\right)\ar`r[u] `u[l] `[lllu] `[ull] [ull]\\
		& & & 0\ar`r[u] `u[l] `[lllu] `[ull] [ull]
	}\quad\xymatrix{
		& 0 & & \\
		& H^{2}\left(N_{0}\right)\ar[r]^{\cup e_{0}} & H^{4}\left(N_{0}\right)\ar@{->}[r]^{\left(p_{0}\right)^{*}} & H^{4}\left(M_{0}\right)\ar`r[u] `u[l] `[lllu] `[ull] [ull] \\
		& 0 & 0\ar[r] & H^{3}\left(M_{0}\right)\ar@{->}`r[u] `u[l] `[lllu] `[ull] [ull] \\
		& H^{0}\left(N_{0}\right)\ar@{->}[r]^{\cup e_{0}} & H^{2}\left(N_{0}\right)\ar@{->}[r]^{\left(p_{0}\right)^{*}} & H^{2}\left(M_{0}\right)\ar`r[u] `u[l] `[lllu] `[ull] [ull]\\
		& & & 0\ar`r[u] `u[l] `[lllu] `[ull] [ull]
	}$$
	Examining the cohomology rings of $M$ and $M_{0}$, we see that $M_{0}$ is a cohomology $S^{2}\times S^{5}\#S^{3}\times S^{4}$. Now $M$ is the connected sum of $M_{0}$ and $(l-1)S^{3}\times S^{4}$. Since $M$ and $S^{3}\times S^{4}$ are both string, we deduce that $M_{0}$ is also string by comparing the characteristic classes. Now $M_{0}$ is a simply connected string cohomology $S^{2}\times S^{5}\#S^{3}\times S^{4}\#\Sigma$ that admits a free circle action. By Propotision \ref{Proposition: homeom S2 times S5 sum S3 tims S4} the manifold $M_{0}$ is diffeomorphic to $S^{2}\times S^{5}\#S^{3}\times S^{4}\#\Sigma'$ for some homotopy $7-$sphere $\Sigma'$ which admits a free circle action. Now we have $S^{2}\times S^{5}\#lS^{3}\times S^{4}\#\Sigma=S^{2}\times S^{5}\#lS^{3}\times S^{4}\#\Sigma'$. Comparing the $\mu-$invariants we see $\mu(\Sigma)=\mu(\Sigma')$. Hence $\Sigma=\Sigma'$ admits a free circle action. Now we conclude that if $S^{2}\times S^{5}\#lS^{3}\times S^{4}\#\Sigma$ admits a free circle action, so does $\Sigma$. This completes the proof of Proposition \ref{Prop: general}, statement 3.
\end{pf}
\begin{pf}[of Lemma \ref{Lemma: k>1, direct construction}]
	For $u\in\mathbb{N}$ we denote $N_{u}=N_{\Sigma}\#u\mathbb{C}P^{3}$, $e_{u}=e_{0}+h_{1}+\cdots+h_{u}\in H^{2}\left(N_{u}\right)$ and $N_{0}=N_{\Sigma}$. Then $N_{u}$ is spin and $e_{k}$ is primitive. In both case $e\in H^{2}(N)$ is primitive, and we identify $e$ with $e_{k-\varepsilon}\in H^{2}\left(N_{k-1-\varepsilon}\right)$ via the isomorphism $H^{2}(N)\cong H^{2}\left(N_{k-1-\varepsilon}\right)$, $\varepsilon=0$, $1$. By \cite[Theorem B]{duan2022circle} if $l$ is even we have $N_{e}=\left(N_{k-1}\right)_{e_{k-1}}\#\Sigma\left(\frac{l}{2}S^{3}\times S^{3}\right)=\left(N_{k-1}\right)_{e_{k-1}}\#l S^{3}\times S^{4}$, and if $l$ is odd $N_{e}$ is diffeomorphic to $\left(N_{k-2}\right)_{e_{k-2}}\#\Sigma\left(S^{2}\times S^{4}\#\frac{l-1}{2}S^{3}\times S^{3}\right)=\left(N_{k-2}\right)_{e_{k-2}}\# S^{2}\times S^{5}\#l S^{3}\times S^{4}$. Now we compute $\left(N_{u}\right)_{e_{u}}$. Since $N_{u}$ is spin, by \cite[Theorems A and B(1)]{galazgarcía2023free} we have $\left(N_{u+1}\right)_{e_{u+1}}=\left(N_{u}\right)_{e_{u}}\#\widetilde{\Sigma}_{h}\mathbb{C}P^{3}=\left(N_{u}\right)_{e_{u}}\# S^{2}\times S^{5}$. Hence by induction $\left(N_{u}\right)_{e_{u}}=(u+1)S^{2}\times S^{5}\#\Sigma$. Now we conclude that $N_{e}=kS^{2}\times S^{5}\#lS^{3}\times S^{4}\#\Sigma$.
\end{pf}
\begin{rmk}
	From the proof we can say more about spinability of the orbit. If $kS^{2}\times S^{5}\#lS^{3}\times S^{4}\#\Sigma$ admits a free circle action, then it admits a free circle action with a spin orbit. While from the proof it is unclear in general whether we can have a nonspin orbit.
\end{rmk}
\begin{pf}[of Theorem \ref{main}]
	 The cases 4 and 5 of Theorem \ref{main} follow from \cite[Theorem 1.3]{RegCirAct2cnt7mfd}. The case 1 of Theorem \ref{main} follows from Proposition \ref{Prop: general}, statement 1. The case 2 of Theorem \ref{main} follows from Proposition \ref{Propotision: homeom S2 times S5} and Proposition \ref{Prop: general}, statement 2. And the case 3 of Theorem \ref{main} follows from Proposition \ref{Proposition: homeom S2 times S5 sum S3 tims S4} and Proposition \ref{Prop: general}, statement 3. Now we complete the proof of Theorem \ref{main}.
\end{pf}

\begin{rmk}
	In Section \ref{More cnt sum} we use the new technique suspension operation to study the cases $k>1$ or $l>1$. We did not follow the approach we used in Sections \ref{Homeom S2 times S5} and \ref{Homeom S2 times S5 sum S3 times S4}. One reason is when the second Betti number becomes greater, the invariants for $6-$manifolds and $7-$manifolds become much more complicated. Another reason is that so far $s-$invariant for $7-$manifolds is merely an invariant of \emph{polarized diffeomorphism}. This means that two such manifolds have the same invariant if and only if there is a diffeomorphism between them that \emph{preserves the prescribed ordered basis of second cohomology groups}. In general this is not a diffeomorphism invariant when the second Betti number is greater than $1$.
\end{rmk}

\section{Miscellaneous complements}\label{Miscellaneous}

If we only apply the method we used in Sections \ref{Homeom S2 times S5} and \ref{Homeom S2 times S5 sum S3 times S4}, we can still deduce some results which are not contained in Theorem \ref{main}. Proposition \ref{Prop: k even, kS2 times S5 only spin orbit} is a conclusion on the orbit of free circle action on spin cohomology $k S^{2}\times S^{5}$, and  Proposition \ref{Prop: autom of k S2 times S5} concerns the self-diffeomorphisms of $k S^{2}\times S^{5}$. 
\begin{prop}\label{Prop: k even, kS2 times S5 only spin orbit}
	Let $k$ be a positive even integer. A spin cohomology $kS^{2}\times S^{5}$ admits NO free circle action with nonspin orbit. %If a spin cohomology $kS^{2}\times S^{5}$ admits a free action with spin orbit $N$ and $e$ is the Euler class of the associated circle bundle, then $\mu(e,e,e)$ must be odd.
\end{prop}
\begin{pf}
	Let $M$ be a spin cohomology $kS^{2}\times S^{5}$ that admits a free circle action. Let $N$ be such an orbit. Let $M\xrightarrow{p}N$ be the bundle projection and let $e\in H^{2}(N)$ be the Euler class of associated circle bundle. From the long exact sequence of homotopy groups and Gysin sequence, we can deduce that $N$ is a simply connected $6-$manifold whose only nontrivial cohomology groups are $H^{0}(N)\cong H^{6}(N)\cong\mathbb{Z}$, $H^{2}(N)\cong H^{4}(N)\cong\mathbb{Z}^{k+1}$ and $H^{2}(N)\xrightarrow{\cup e}H^{4}(N)$ is an isomorphism. Extend $e_{0}=e$ to an ordered basis $\left(e_{0},e_{1},\cdots,e_{k}\right)$ of $H^{2}(N)$. Equip $H^{4}(N)$ with the dual basis $\left(e^{0},\cdots,e^{k}\right)$. The matrix representation of $H^{2}(N)\xrightarrow{\cup e}H^{4}(N)$ with respect to the given bases is $M_{e}=\left(\mu\left(e_{0},e_{i},e_{j}\right)\right)_{0\leqslant i,j\leqslant k}$, and $M_{e}$ is symmetric and invertible. The second Stiefel-Whitney class of $N$ is computed as before. From the relation $w_{2}(M)=p^{*}\left(w_{2}(N)+\overline{e}\right)=0$ and $\mathrm{mod}$ $2$ Gysin sequence we deduce that $w_{2}(N)=0$ or $\overline{e}$.
	
	Assume $N$ is nonspin, $w_{2}(N)=\overline{e}$. Formula \eqref{Formula: Wu formula} implies that $\mathrm{Sq}^{2}\left(\overline{e_{i}}^{2}\right)=\overline{e_{0}}\cdot\overline{e_{i}}^{2}$, $0\leqslant i\leqslant k$. Hence $\overline{\mu\left(e_{0},e_{i},e_{i}\right)}=0\in\mathbb{Z}/2$, $\overline{M_{e}}\in\mathrm{GL}(k+1,\mathbb{Z}/2)$ has vanishing diagonal and $\overline{M_{e}}$ admits a skew-symmetric integral lifting $A\in\mathrm{GL}(k+1,\mathbb{Z})$. By assumption $k$ is even, thus $k+1$ is odd. Hence $\mathrm{det}A=0$ and $\overline{\mathrm{det}M_{e}}=\mathrm{det}\overline{M_{e}}=\mathrm{det}\overline{A}=\overline{\mathrm{det}A}=0\in\mathbb{Z}/2$. While $M_{e}$ is invertible, $\mathrm{det}M_{e}=\pm1$ and $\overline{\mathrm{det}M_{e}}=\overline{1}\in\mathbb{Z}/2$. This is a contradiction. Hence a spin cohomology $kS^{2}\times S^{5}$ cannot admit free circle action with nonspin orbit.
\end{pf}
\begin{prop}\label{Prop: autom of k S2 times S5}
	For any $g\in\mathrm{GL}(k,\mathbb{Z})$ there is a self-diffeomorphism $\varphi$ of $k S^{2}\times S^{5}$ such that $g=H^{2}(\varphi)$.
\end{prop}
\begin{pf}
	We note that $M=kS^{2}\times S^{5}$ satisfies the condition of \cite[Theorem 6]{SurgeryAndDuality}. Let $W$ be the boundary connected sum of $k$ copies of $S^{2}\times D^{6}$. Then is a spin coboundary of $M$ such that boundary inclusion $M\xrightarrow{\iota}W$ induces an isomorphism on the second integral cohomology groups and $H^{4}(W,M;\mathbb{Q})=0$. In particular $\sigma(W,M)=0$ and for any two elements $z,w\in H^{4}(W,M;\mathbb{Q})$, $\left<z\cup w,[W,M]\right>=0$. Let $x=\left(x_{1},\cdots,x_{k}\right)$, $y=\left(y_{1},\cdots,y_{k}\right)$ be two ordered bases of $H^{2}(M)$ and let $g$ be the transition matrix from $x$ to $y$. Since $H^{2}(W)\xrightarrow{\iota^{*}}H^{2}(M)$ is an isomorphism, we can take $\widetilde{x_{i}}$, $\widetilde{y_{i}}\in H^{2}(W)$ such that $\iota^{*}\widetilde{x_{i}}=x_{i}$, $\iota^{*}\widetilde{y_{i}}=y_{i}$. Set $\widetilde{x}=\left(\widetilde{x_{1}},\cdots,\widetilde{x_{k}}\right)$, $\widetilde{y}=\left(\widetilde{y_{1}},\cdots,\widetilde{y_{k}}\right)$. Then $\left(W;\widetilde{x}\right)\sqcup\left(-W;\widetilde{y}\right)$ is a spin bordism between $(M;x)$ and $(M;y)$, and its characteristic numbers that are required to vanish in \cite[Theorem 6]{SurgeryAndDuality} are zero. Here the oriented manifold $-W$ and $W$ have the same underlying manifold and opposite orientations. Hence $(M;x)$ and $(M;y)$ are polarized diffeomorphic, i.e. there is a self-diffeomorphism $\varphi$ of $M$ such that $\varphi^{*}x=y$, so that $\varphi^{*}=g$.
\end{pf}
\begin{rmk}
	\cite[Corollary 4.7]{galazgarcía2023free} implies this particular spin case in dimension $7$.
\end{rmk}
\begin{ackn}
	The author would like to express gratitude to Yi Jiang for introducing the topic and for her insightful discussions and suggestions. Thanks are also due to Yang Su for his valuable suggestions, and to Haibao Duan for his thoughtful comments. The author's research was partially supported by NSFC 11801298.
\end{ackn}

%\section*{Acknowledgement}
%The author would thank my supervisor professor Yi Jiang. She offered this problem to me, patiently listened to my report, gave me precious suggestions and enlightened my road of research, The author would also thank professor Yang Su and professor Weiyan Chen for useful comments. Finally, the author would thank Chao Shen, Hao Liang, Liyuan Ye and Xurui Zhao for helpful discussion on those algebraic and number-theoretic problems that occurred during the study.

\bibliography{FreeCircleActionsOnCertainSimplyConnected7Manifolds}

\begin{thebibliography}{GGR23}

\bibitem[Bru69]{Brumfiel69}
G.~Brumfiel.
\newblock Differentiable ${S}^{1}$ actions on homotopy spheres, 1969.

\bibitem[DL05]{Duan2005CircleBO}
Haibao Duan and Chao Liang.
\newblock Circle bundles over 4-manifolds.
\newblock {\em Archiv der Mathematik}, 85:278--282, 2005.

\bibitem[Dua22]{duan2022circle}
Haibao Duan.
\newblock Circle actions and suspension operations on smooth manifolds, eprint
  arXiv: 2202.06225, 2022.

\bibitem[EK62]{EKinv}
James Eells and Nicolaas~H. Kuiper.
\newblock An invariant for certain smooth manifolds.
\newblock {\em Annali di Matematica Pura ed Applicata}, 60:93--110, 1962.

\bibitem[Gei08]{Geiges_2008}
Hansjörg Geiges.
\newblock {\em An Introduction to Contact Topology}.
\newblock Cambridge Studies in Advanced Mathematics. Cambridge University
  Press, 2008.

\bibitem[GGR23]{galazgarcía2023free}
Fernando Galaz-García and Philipp Reiser.
\newblock Free torus actions and twisted suspensions, eprint arXiv: 2305.06068,
  2023.

\bibitem[Jia14]{RegCirAct2cnt7mfd}
Yi~Jiang.
\newblock {Regular circle actions on 2-connected 7-manifolds}.
\newblock {\em Journal of the London Mathematical Society}, 90(2):373--387, 07
  2014.

\bibitem[Jup73]{Jupp1973}
P.~E. Jupp.
\newblock Classification of certain 6-manifolds.
\newblock {\em Mathematical Proceedings of the Cambridge Philosophical
  Society}, 73(2):293–300, 1973.

\bibitem[KM63]{KervaireMilnor63}
Michel~A. Kervaire and John~W. Milnor.
\newblock Groups of homotopy spheres: I.
\newblock {\em Annals of Mathematics}, 77(3):504 -- 537, 1963.

\bibitem[Kre99]{SurgeryAndDuality}
Matthias Kreck.
\newblock Surgery and duality.
\newblock {\em Annals of Mathematics}, 149(3):707--754, 1999.

\bibitem[Kre18]{Kreck2018OnTC}
Matthias Kreck.
\newblock On the classification of 1‐connected 7‐manifolds with torsion
  free second homology.
\newblock {\em Journal of Topology}, 11, 2018.

\bibitem[KS88]{KS88}
Matthias Kreck and Stephan Stolz.
\newblock A diffeomorphism classification of 7-dimensional homogeneous einstein
  manifolds with {SU}(3)×{SU}(2)×{U}(1)-symmetry.
\newblock {\em Annals of Mathematics}, 127(2):373--388, 1988.

\bibitem[KS91]{Kreck1991SomeNH}
Matthias Kreck and Stephan Stolz.
\newblock Some nondiffeomorphic homeomorphic homogeneous 7-manifolds with
  positive sectional curvature.
\newblock {\em Journal of Differential Geometry}, 33:465--486, 1991.

\bibitem[KS98]{Kreck1991SomeNHCorrection}
Matthias Kreck and Stephan Stolz.
\newblock A correction on: ``{S}ome nondiffeomorphic homeomorphic homogeneous
  7-manifolds with positive sectional curvature''.
\newblock {\em Journal of Differential Geometry}, 49(1):203 -- 204, 1998.

\bibitem[Lee12]{JohnLee}
John~M. Lee.
\newblock {\em Introduction to Smooth Manifolds}.
\newblock Graduate Texts in Mathematics. Springer New York, NY, 2012.

\bibitem[MS74]{MilnorStasheff}
John~W. Milnor and James~D. Stasheff.
\newblock {\em Characteristic Classes}.
\newblock Princeton University Press, 1974.

\bibitem[MY68]{MontgomeryYang68}
Deane Montgomery and C.~T. Yang.
\newblock Differentiable actions on homotopy seven spheres, ii.
\newblock In Paul~S. Mostert, editor, {\em Proceedings of the Conference on
  Transformation Groups}, pages 125--134, Berlin, Heidelberg, 1968. Springer
  Berlin Heidelberg.

\bibitem[TW15]{Tuschmann2015ModuliSO}
Wilderich Tuschmann and David~J. Wraith.
\newblock {\em Moduli Spaces of Riemannian Metrics}.
\newblock Oberwolfach Seminars. Birkhäuser Basel, 2015.

\bibitem[Wal66]{Wall66}
C.~T.~C. Wall.
\newblock Classification problems in differential topology. {V}: {On} certain
  6-manifolds.
\newblock {\em Invent. Math.}, 1:355--374, 1966.

\end{thebibliography}

\textsc{Department of Mathematics, Tsinghua University, Beijing, P.R.China}

Email address: \nolinkurl{xufupeng16@mails.ucas.ac.cn}

\end{document}